\documentclass[oneside,a4paper,final]{amsart}
\pdfoutput=1

\usepackage[UKenglish]{babel}
\usepackage[utf8]{inputenc}
\usepackage{lmodern}
\usepackage{amssymb,amsfonts}
\usepackage{amsthm}
\usepackage{amsmath}
\usepackage{amsmath}
\usepackage{enumerate}
\usepackage{mathtools}
\usepackage{color}
\usepackage{tikz}
\usetikzlibrary{decorations.pathreplacing}
\usetikzlibrary{arrows}
 \usepackage{graphicx}
\usepackage{ stmaryrd }
\usetikzlibrary{hobby}
\usepackage{dsfont}
\usepackage[backend=bibtex,firstinits,maxnames=4,style=alphabetic,isbn=false, date=year, doi=false, eprint= true, url= false]{biblatex}
\usepackage{esint}

\usepackage{microtype}

\usepackage[notcite,notref,color]{showkeys}
\definecolor{labelkey}{gray}{.6}
\definecolor{refkey}{gray}{.6}

\usepackage[colorlinks=true, pdfstartview=FitV, linkcolor=blue, citecolor=blue, urlcolor=blue,pagebackref=false]{hyperref}

\usepackage[margin=1.3in,headheight= 13.6 pt, asymmetric]{geometry}

\newcommand{\Bcal}{\mathcal{B}}

\newcommand{\Dcal}{\mathcal{D}}

\newcommand{\Rcal}{\mathcal{R}}

  \providecommand{\eta}{{\boldsymbol{\eta}}}

\DeclareMathOperator{\dist}{dist}

\newcommand{\norm}[1]{\|#1\|}

\newcommand{\normB}[1]{\Bigl\|#1\Bigr\|}

\newcommand{\abs}[1]{|#1|}

\newcommand{\absB}[1]{\Bigl|#1\Bigr|}

\newcommand{\dpr}[1]{\langle #1 \rangle}

\newcommand{\N}{\mathbb{N}}
\newcommand{\R}{\mathbb{R}}

\newcommand{\eps}{\varepsilon}
\newcommand{\lrangle}[1]{\langle #1 \rangle}

\author{Richard  M. H\"ofer}
\author{Karina  Kowalczyk}
\author{Sebastian Schwarzacher}

\title[Low Mach and homogenization limit of compressible fluids]{Darcy's law as low Mach and homogenization limit of a compressible fluid in perforated domains}

\date{\today}

\theoremstyle{definition}
\newtheorem{defi}{Definition}[section]
\theoremstyle{plain}
\newtheorem{theorem}[defi]{Theorem}
\newtheorem{conj}[defi]{Conjecture}
\newtheorem{lem}[defi]{Lemma}
\newtheorem{prop}[defi]{Proposition}

\theoremstyle{definition}
\newtheorem{rem}[defi]{Remark}

\let\div\undefined
\DeclareMathOperator{\div}{div}
\DeclareMathOperator{\supp}{supp}

\newcommand{\rhoeps}{\rho_{\varepsilon}}
\newcommand{\ueps}{u_{\varepsilon}}
\newcommand{\peps}{p_{\varepsilon}}

\newcommand{\sigmaeps}{\sigma_{\varepsilon}}
\newcommand{\omegaeps}{\Omega_\varepsilon}
\newcommand{\rhotilde}{\tilde{\rho}_{\varepsilon}}

\newcommand{\utilde}{\tilde{u}_{\varepsilon}}

\newcommand{\weak}{\rightharpoonup}

 \newcommand{\wk}{w_k^{\varepsilon}}
\newcommand{\qk}{q^{\varepsilon}_k}

\newcommand{\inteps}{\int_{\omegaeps}}

\numberwithin{equation}{section} 

\newcommand{\e}{{\varepsilon}}
\newcommand{\weakto}{\rightharpoonup}

\mathtoolsset{showonlyrefs}

\bibliography{bibl_master_thesis}

\begin{document}

\begin{abstract}
We consider  the homogenization limit of the compressible barotropic Navier-Stokes equations in a three-dimensional domain perforated by periodically distributed identical particles.
We study the regime of particle sizes and distances such that the volume fraction of particles tends to zero but their resistance density tends to infinity.
Assuming that the Mach number is decreasing with a certain rate,  the rescaled velocity and pressure of the microscopic system converges to the solution of an effective equation which is given by Darcy's law. The range of sizes of particles we consider are exactly the same which lead to Darcy's law in the homogenization limit of incompressible fluids.
Unlike  previous results for the Darcy regime we estimate the deficit related to the pressure approximation via the \texorpdfstring{Bogovski\u{ı}}{Bogovski} operator.
This allows for more flexible estimates of the pressure in Lebesgue and Sobolev spaces and allows to proof convergence results for all barotropic exponents $\gamma> \frac32$. 
\end{abstract}

\maketitle




\section{Introduction}

\subsection{Setting and description of the main result}

\begin{figure}[h] 
\begin{center}

\begin{tikzpicture} [scale = 0.4]
 \filldraw[fill={rgb:black,1;white,2}] (0,0) to [closed, curve through = {(1.2,2) (0.4,3.3) (1.9,4.8)  (3.9,4.2)(3.2,2.4) (4,0) (3,-1.4) (0.3,-1.2)}] (0,0);
 \node (A) at (2.2:1.5) [above right] {$T$};
 
 \node (B) at (-2.5, 4.0) [above right] {$Q$};
 
 \draw[draw=black] (-2.8,-3.5) rectangle ++(10,10);
 
\end{tikzpicture}
\end{center}
\caption{A single cell}  \label{fig:single.particle}
\end{figure}
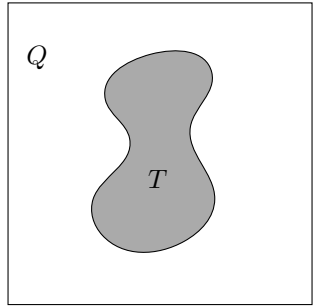

Let $\Omega$ be a bounded connected open set in $\mathbb{R}^3$, with smooth boundary $\partial \Omega$, $\Omega$ being locally located  on one side of its boundary. Let $\varepsilon$ be a sequence of strictly positive real numbers tending to zero. Let $Q = (-1,1)^3$ and $B_1 \subset Q$ the unit ball. Let $T \Subset B_1$, the reference particle, be a fixed closed set with smooth boundary, such that $B_1 \setminus T$ is a connected open set, locally located on one side of its boundary (see Figure \ref{fig:single.particle}).

We cover the set $\Omega$ with a regular mesh of size $2 \varepsilon$ and we denote by $x_i^{\varepsilon}$ the center of the cell with index $i$ at level $\varepsilon$. Let $P_i^{\varepsilon} = x_i^{\varepsilon} + (- \varepsilon, + \varepsilon)^3$ be the cell with center $x_i^{\varepsilon}$. Let $ i \in \{ 1, \dots, N(\varepsilon)\}$ be those indices for which the cell $P_i^{\varepsilon}$ is entirely included in $\Omega$. Let $\alpha \geq 1$ and consider in each of the cells the particle  \begin{equation}
T_i^{\varepsilon} := x_i^{\varepsilon} + \varepsilon^{\alpha} T , \hspace{2 ex} i = 1 \dots N(\varepsilon).
\end{equation}
Now we define the perforated domain $\omegaeps$ by \begin{equation}\label{omegaeps}
\omegaeps = \Omega \setminus \bigcup_{ i = 1}^{N(\varepsilon)} T_i^{\varepsilon}.
\end{equation}

All the objects defined above are illustrated in Figure \ref{fig:perforated.domain}.

\begin{figure} 
\hspace{ 2.5 cm}
\begin{tikzpicture}[scale = 1.3 ]
\filldraw[fill=white, thick] (0,1.9) to [closed, curve through = {(1.15,1.85) (1.9,1) (1.2,0)  (1.7,-1.3)(0.9,-2.1) (-1,-2) (-2.3,-0.5) (-2,0.5)(-2.2, 1.8) (-1.5, 2.3)(-0.5, 1.9)}] (0,1.9);

\draw [step=0.9,black, ] (-2.5,-2.6) grid (2.2,2.5);

  \node at (-1.35, -0.45){
   	\begin{tikzpicture}[scale = 0.04]
    	\filldraw[fill={rgb:black,1;white,2}] (-2.2,-1.5) to [closed, curve 					through = {(-1,0.5) (-1.8,1.8) (-0.3,3.3)  (1.7,2.7)(1,0.9) 						(1.8,-1.5) (0.8,-2.9) (-1.9,-2.7)}] (-2.2,-1.5);
    \end{tikzpicture}
  };
  
  \node at (-1.35, 0.45){
   	\begin{tikzpicture}[scale = 0.04]
    	\filldraw[fill={rgb:black,1;white,2}] (-2.2,-1.5) to [closed, curve 					through = {(-1,0.5) (-1.8,1.8) (-0.3,3.3)  (1.7,2.7)(1,0.9) 						(1.8,-1.5) (0.8,-2.9) (-1.9,-2.7)}] (-2.2,-1.5);
    \end{tikzpicture}
  };
  
  \node at (-1.35, 1.35){
   	\begin{tikzpicture}[scale = 0.04]
    	\filldraw[fill={rgb:black,1;white,2}] (-2.2,-1.5) to [closed, curve 					through = {(-1,0.5) (-1.8,1.8) (-0.3,3.3)  (1.7,2.7)(1,0.9) 						(1.8,-1.5) (0.8,-2.9) (-1.9,-2.7)}] (-2.2,-1.5);
    \end{tikzpicture}
  };

  \node at (-0.45, -1.35){
   	\begin{tikzpicture}[scale = 0.04]
    	\filldraw[fill={rgb:black,1;white,2}] (-2.2,-1.5) to [closed, curve 					through = {(-1,0.5) (-1.8,1.8) (-0.3,3.3)  (1.7,2.7)(1,0.9) 						(1.8,-1.5) (0.8,-2.9) (-1.9,-2.7)}] (-2.2,-1.5);
    \end{tikzpicture}
  };
  
  \node at (0.45, -1.35){
   	\begin{tikzpicture}[scale = 0.04]
    	\filldraw[fill={rgb:black,1;white,2}] (-2.2,-1.5) to [closed, curve 					through = {(-1,0.5) (-1.8,1.8) (-0.3,3.3)  (1.7,2.7)(1,0.9) 						(1.8,-1.5) (0.8,-2.9) (-1.9,-2.7)}] (-2.2,-1.5);
    \end{tikzpicture}
  };

  \node at (-0.45, -0.45){
   	\begin{tikzpicture}[scale = 0.04]
    	\filldraw[fill={rgb:black,1;white,2}] (-2.2,-1.5) to [closed, curve 					through = {(-1,0.5) (-1.8,1.8) (-0.3,3.3)  (1.7,2.7)(1,0.9) 						(1.8,-1.5) (0.8,-2.9) (-1.9,-2.7)}] (-2.2,-1.5);
    \end{tikzpicture}
  };
  
  \node at (-0.45, 0.45){
   	\begin{tikzpicture}[scale = 0.04]
    	\filldraw[fill={rgb:black,1;white,2}] (-2.2,-1.5) to [closed, curve 					through = {(-1,0.5) (-1.8,1.8) (-0.3,3.3)  (1.7,2.7)(1,0.9) 						(1.8,-1.5) (0.8,-2.9) (-1.9,-2.7)}] (-2.2,-1.5);
    \end{tikzpicture}
  };

 \node (C) at (-2.45,2.1) {$\omegaeps$};
 
  \node at (0.45, -0.45){
   	\begin{tikzpicture}[scale = 0.04]
    	\filldraw[fill={rgb:black,1;white,2}] (-2.2,-1.5) to [closed, curve 					through = {(-1,0.5) (-1.8,1.8) (-0.3,3.3)  (1.7,2.7)(1,0.9) 						(1.8,-1.5) (0.8,-2.9) (-1.9,-2.7)}] (-2.2,-1.5);
    \end{tikzpicture}
  };
  \draw (0.45, -0.45) -- (2.2, -0.45);
  \node (K) at (3.3,-0.45) {$T_i^{\varepsilon} = x_i^{\varepsilon} +\varepsilon^{\alpha} T$};
  \draw (0.65, -0.7) -- (2.4, -1.2);
  \node (J) at (2.7,-1.2) {$P_i^{\varepsilon}$};
  
  \node at (0.45, 0.45){
   	\begin{tikzpicture}[scale = 0.04]
    	\filldraw[fill={rgb:black,1;white,2}] (-2.2,-1.5) to [closed, curve 					through = {(-1,0.5) (-1.8,1.8) (-0.3,3.3)  (1.7,2.7)(1,0.9) 						(1.8,-1.5) (0.8,-2.9) (-1.9,-2.7)}] (-2.2,-1.5);
    \end{tikzpicture}
  };
  
  \node at (0.45, 1.35){
   	\begin{tikzpicture}[scale = 0.04]
    	\filldraw[fill={rgb:black,1;white,2}] (-2.2,-1.5) to [closed, curve 					through = {(-1,0.5) (-1.8,1.8) (-0.3,3.3)  (1.7,2.7)(1,0.9) 						(1.8,-1.5) (0.8,-2.9) (-1.9,-2.7)}] (-2.2,-1.5);
    \end{tikzpicture}
  };
  
  \node at (-0.45, 1.35){
   	\begin{tikzpicture}[scale = 0.04]
    	\filldraw[fill={rgb:black,1;white,2}] (-2.2,-1.5) to [closed, curve 					through = {(-1,0.5) (-1.8,1.8) (-0.3,3.3)  (1.7,2.7)(1,0.9) 						(1.8,-1.5) (0.8,-2.9) (-1.9,-2.7)}] (-2.2,-1.5);
    \end{tikzpicture}
  };
  
  \node (S) at (2.3, 1){};
   \node (EN) at (2.3, -0.1){};
  
\draw[thick,decorate,decoration={brace,amplitude=3pt}]
(S) -- (EN) node[midway, right=4pt, text width=80pt] {$ 2\varepsilon$};

\end{tikzpicture}
\caption{The perforated domain $\Omega_\eps$} \label{fig:perforated.domain}
\end{figure}



We study the homogenization problem for a compressible viscous fluid at low Mach number in the the perforated domain $\omegaeps$. Assuming a  barometric pressure law $P(\rho) = \rho^\gamma$, $\gamma > 0$,  we model the fluid by the compressible Navier-Stokes system
\begin{equation} \label{goal}
\begin{cases}
\sigmaeps^2 \partial_t(\rhoeps \ueps) +\div(\rhoeps \ueps \otimes \ueps) - \Delta \ueps  - \nabla \div \ueps + \frac{1}{\eps^\beta}\nabla\rhoeps^{\gamma} = \rhoeps f + g \quad & \text{in } \omegaeps \times (0,T),\\
\sigmaeps^2 \partial_t \rhoeps +\div(\rhoeps \ueps) = 0 \quad & \text{in } \omegaeps \times (0,T,)\\
\ueps = 0 \quad & \text{on } \partial\omegaeps \times (0,T),\\
\rhoeps  = \rho_{\varepsilon0}, \hspace{1 ex} \rhoeps \ueps  = q_{\varepsilon 0} \quad & \text{on } \Omega_\eps \times \{0\},
\end{cases}
\end{equation}
where $\beta > 0$, and $f,g,\rho_{\eps0}, q_{\eps0}$ are given source terms and initial data. Moreover,
\begin{equation} \label{sigma}
	\sigma_\eps = \e^{\frac{3-\alpha}2},
\end{equation}
which has turned out to be a decisive quantity for such homogenization problems.
 We study the homogenization limit $\eps \to 0$ for $\alpha \in (1,3)$,  $\gamma > 3/2$, and all $\beta$ sufficiently large (depending on $\alpha$ and $\gamma$). We prove convergence for the rescaled fluid velocity $\sigma_\eps^{-2} u_\eps$ and  pressure $p_\eps = \eps^{-\beta} (\rho^\gamma - \lrangle {\rho_\eps^\gamma}_\eps)$ to  limiting functions $(u,p)$. Here $\lrangle {\cdot}_\eps$ denotes the average over $\omegaeps$.
We show that $(u,p)$  is a solution to the following Darcy's law
\begin{equation}\label{darcy_law1}
\begin{cases} u =  \Rcal^{-1}(\rho_0 f + g - \nabla p) \quad & \text{in } \Omega \times (0,T), \\
\div u =  0 \hspace{1 ex} \quad & \text{in } \Omega \times (0,T) ,\\
u\cdot n = 0 \hspace{1 ex} \quad & \text{on } \partial \Omega \times (0,T),
\end{cases}
\end{equation} 
where the resistance matrix $\Rcal$ will be introduced in the following subsection.
 We remark that as $ u_\eps \sim \sigma_\eps^2$, the appearance of the term $\sigma_\eps^2$ in front of the time derivatives is due to the natural time scale of the system.

We postpone the precise statement of our main result to Section \ref{sec:main.results}.
The first main result (Theorem~\ref{theorem1}) is for the steady compressible Stokes system where we proof the convergence to the steady Darcy's law for {\em all barotropic exponents} $\gamma\geq 1$. The second main result (Theorem~\ref{theorem3}) provides the analogous convergence result for the unsteady compressible Navier-Stokes system for exponents $\gamma>\frac{3}{2}$.

It seems to us that since the seminal work of Massmoudi~\cite{Masmoudi}, homogenization results for compressible fluids in the
Darcy regime are rather sparse. In particular we adress here for the first time the regime $\alpha \in (1,3)$.  To our knowledge, another novelty of our result 
among compressible homogenization problems is the treatment of all barotropic exponents  $\gamma>\frac{3}{2}$. As is apparent from Theorem~\ref{theorem1}, the reason of the restrictions on the barotropic exponents is a question of (missing) regularity estimates. 
In this paper, we succeed to deal with all exponents $\gamma>\frac32$ due to a noval treatment of the pressure approximation $p_\eps$. More precisely, we decompose this pressure through the dual of the so called Bogovski\u{\i} operator~\cite{bogovskij} into components corresponding to the inertial forces and to the viscous and external forces. As such, the decomposition inherits the (recently developed) estimates in the full range of exponents for Lebesgue and Sobolev spaces in perforated domains. These features might well be of independent interest for future applications.

\subsection{Heuristics leading to the Darcy's law} \label{sec:heuristics}

We analyze the homogenization problem $\eps \to 0$ in the regime, determined by the parameters $\alpha, \beta$ and $\gamma$, where, locally at each particle, the fluid flow is well approximated by the incompressible steady Stokes equations. 
Then, the homogeneous Dirichlet boundary conditions at the particles give rise to a
so-called Brinkman force term $\sigma_\eps^{-2} \bar \Rcal u_\eps$ that is determined by the resistance density $\sigma_\eps^{-2} \bar \Rcal$ of the particles, where $\bar \Rcal$ is the resistance matrix of the reference particle $T$.

More precisely,  since the particles are far away from each other, at a fixed time $t \in (0,T)$,
$u_\eps$ around particle $T_i^\eps$ can be approximated by the solution $v$ to
\begin{align}
	\begin{cases}
		- \Delta v + \nabla q = 0 \quad & \text{in } \R^3 \setminus T_i^\eps, \\
		\div v = 0 \quad & \text{in } \R^3 \setminus T_i^\eps, \\
		v = 0 \quad  & \text{on } \partial T_i^\eps, \\
		\lim_{|x - x_i^\eps| \to \infty} v(x)   = \bar u_\eps(x_i^\eps),
	\end{cases}
\end{align}
where $\bar u_\eps(x_i^\eps)$ is the (average) value of $u_\eps$ \emph{around} particle $T_i^\eps$, i.e. at distances $\eps^\alpha \ll |x - x_i^\eps|  \ll \eps$.

There is a linear relation between $u_\eps(x_i^\eps)$ and the drag force
\begin{align}
F = \int_{\partial T_i^\eps} (\nabla v + (\nabla v)^T) n - q n,
\end{align}
namely
\begin{align} \label{eq:force.single.particle}
	F = \eps^\alpha \bar{\Rcal} (\bar u_\eps(x_i^\eps)).
\end{align}
This is immediate from scaling considerations and the definition of $\bar \Rcal$ via
\begin{align} \label{def.Rcal.bar}
	\bar \Rcal_{jk} =  e_j  \cdot \int_{\partial T_i^\eps} (\nabla w_k + (\nabla w_k)^T) n - q_k n,
\end{align}
 where $(w_k,q_k)$ is the solution to the Stokes problem
 \begin{align} \label{stokessingle}
\begin{cases}
\nabla q_k - \Delta w_k= 0  \quad  & \text{in } \mathbb{R}^N\setminus T,\\
\div w_k = 0  \quad & \text{in } \mathbb{R}^N\setminus T,\\
w_k = 0 \quad  & \text{on } \partial T,\\
w_k = e_k \quad  & \text{at infinity.}
\end{cases}
\end{align}
In view of \eqref{eq:force.single.particle} the effect of the particles is
approximated by replacing the Dirichlet boundary data by an additional force $-(2\eps)^{-3} \eps^{\alpha} \bar \Rcal u_\eps = - \frac{1}{8}\sigma_\eps^{-2} \bar \Rcal u_\eps $, where the factor $(2\eps)^{-3}$ reflects the particle density. This leads to the approximate momentum equation 
\begin{align}
	\sigmaeps^2 \partial_t(\rhoeps \ueps) +\div(\rhoeps \ueps \otimes \ueps) - \Delta \ueps  - \nabla \div \ueps  + \frac{\bar{\Rcal} u_\eps}{8 \sigma_\eps^2} + \frac{1}{\eps^\beta}\nabla\rhoeps^{\gamma} = \rhoeps f + g.
\end{align}
Thus, $u_\eps \sim \sigma_\eps^2$. Moreover,  $\sigma_\eps \to 0$ since  we assume $\alpha < 3$, which entails that the Brinkman force is much larger than all other terms involving $u_\eps$. Combining this with the low Mach assumption as $\eps \to 0$,
we formally arrive at the Darcy's law \eqref{darcy_law1}.
Here we denote
\begin{align} \label{def:Rcal}
	\Rcal = \frac{1}{8} \bar \Rcal
\end{align}
to avoid the appearance of the factor $1/8$.

\subsection{Previous results}

There is a vast literature on homogenization of viscous fluid flows in perforated domains with Dirichlet boundary conditions.
We will not provide an exhaustive discussion of this literature here. More references can be found in the cited papers.

The results can be broadly characterized by two aspects, the fluid equations that are considered and the assumptions on the particle configuration.
On the one hand, the fluid equations range from stationary incompressible Stokes equations to the incompressible and compressible Navier-Stokes equations. 
On the other hand, the particle configurations range from periodic arrangements as considered here to more general and even stochastic particle distributions.

\smallskip

The most accomplished theory is available for the incompressible steady Stokes and Navier-Stokes equations. To our knowledge, the first rigorous result is due to Tartar
\cite{tartar}, who established a Darcy's law for incompressible steady Stokes equations for $\alpha = 1$. We emphasize that the case $\alpha = 1$ is fundamentally different from $\alpha > 1$. This is due to the fact that $\eps$ is the only small length scale in this case which makes the problem more accesible but also means that the interaction between the particles is significant in this regime.
Therefore, the resistance matrix $\Rcal$ in the Darcy's law \eqref{darcy_law1} needs to be replaced by a different resistance matrix $ \Rcal_p$ in this case, which is determined by solving a Stokes problem in a cell with periodic boundary conditions.

For the incompressible Stokes equations, all the cases $\alpha > 1$ have been addressed by Allaire in 
\cite{allaire1} and \cite{allaire2} for periodically distributed particles.
The results for the Stokes equations are summarized in Table \ref{tab:Allaire}. By compactness, they generalize to the steady Navier-Stokes equations.
\begin{table}[h]
\begin{center}
\begin{tabular}{c|c|c}
small particles & critical size of particles & large particles\\[0.3 cm]
$\alpha >3$ & $\alpha =3$ & $1 < \alpha < 3 $ \\
\hline 
 & & \\
 $- \Delta u + \nabla p = f$ & $- \Delta u + \nabla p + \Rcal u = f$ & $u = \Rcal^{-1}(f - \nabla p)$\\
$\div(u) = 0$ & $\div(u) = 0$ & $\div(u) = 0$\\ 
 & & \\
\hline
 & & \\
Stokes' law & Brinkmann's law & Darcy's law\\
\end{tabular}
\end{center}
\caption{Homogenization results for incompressible Stokes equations} \label{tab:Allaire} 
\end{table}

In the critical regime $\alpha = 3$, these results have been generalized to non-homogeneous Dirichlet boundary conditions and more general particle configurations,
including stochastic configurations (also with random sizes of the particles) and the study of the fluctuations around the homogenization limit. Mentioning some of these results, one can refer to \cite{Rubinstein86,desvillettes,
GiuntiHoefer19,CarrapatosoHillairet20,HoeferJansen20}
. Very recently, the results by Allaire in the regime $\alpha \in (1,3)$ have been generalized to stochastic particle configurations in \cite{Giunti21}.

In the case $\alpha =1$, Beliaev and Kozlov established a Darcy's law  for the homogenization of the Stokes equations in a randomly perforated domain in \cite{BeliaevKozlov96}.

\smallskip

For $\alpha \in \{1,3\}$, there are also results for the incompressible instationary Navier-Stokes equations.
For $\alpha = 1$ and periodic particle configurations, Mikeli\'c \cite{Mikelic91} established a Darcy's law.  Feireisl, Namlyeyeva and Ne\v{c}asova \cite{FeireislNamlyeyevaNecasova16} proved the appearance of a Brinkman term in the homogenized equation in the critical regime $\alpha = 3$ under rather general assumptions on the particle configurations.

\smallskip

While the homogenization problem for incompressible fluids is quite well understood, there are only few results for compressible fluids. These results are restricted to periodic particle configurations and the cases of $\alpha = 1$ on the one hand, and of very small particles, i.e. $\alpha > 3$, on the other hand.
In the regime of very small particles, the effect of the particles have been shown to be negligible for the stationary and evolutionary Navier-Stokes and Navier-Stokes-Fourier system by Feireisl and Lu \cite{FeireislLu15}, Diening and the former~\cite{DieFeiLu17}, Lu and the third author \cite{schwarzacher}, and Lu and Pokorn\'y \cite{LuPokorny20}.

For $\alpha = 1$, Masmoudi~\cite{Masmoudi} studied the homogenization problem for the compressible 
Navier-Stokes system \eqref{goal} for $\beta = 0$ and $\gamma \geq 3$. He proved that the homogenization limit is given by the compressible version of the Darcy's law, which takes the form of a porous medium equation
\begin{align}
\begin{cases}
\theta\partial_t \rho + \div \left( \rho u ) \right)=  0 & \text{ in } \Omega \times (0,T),\\
u =  \Rcal_p^{-1} (\rho f + g - \nabla \rho^{\gamma}) \text{ on } \{ \rho > 0 \} & \text{ in }\Omega \times (0,T), \\
u \cdot n = 0 &\text{ on } \partial \Omega \times (0,T),
\end{cases}
\end{align}
where $\theta$ is the volume fraction of the particles. This result has been generalized to the Navier-Stokes-Fourier system by Feireisl, Novotn\'y and Takahashi \cite{FeireislNovotnyTakahashi10}.
Again, the key observation for these results is that locally around each particle, the fluid can be approximately treated as an incompressible creeping flow, which leads to the same resistance matrix $ \Rcal_p$ as in the Darcy's law obtained in \cite{tartar} for the Stokes equations. To our knowledge these are the only homogenization results for compressible fluids in the Darcy regime.

\subsection{Significance and open problems}

To our knowledge, the homogenization problem for the compressible  Navier-Stokes equations has not been addressed so far in the regime $1 < \alpha < 3$.
 In view of Allaire's results displayed in Table \ref{tab:Allaire} above, one could naively anticipate that Masmoudi's results from~\cite{Masmoudi} can be generalized appropriately to the regime $1 < \alpha < 3$. However, as we pointed out before, the key ingredient for such a homogenization result is that the fluid near each particle is well approximated by the incompressible Stokes equations. Formally, the local Reynolds number is easily seen to be very small since both the local length scale $\eps^{\alpha}$ and the fluid velocity $\sigma_\eps^2$  are very small. Thus, the convective terms are negligible.
However, for local incompressibility, one needs that the square of the local Mach number divided by the local Reynolds number tends to zero (For the non-dimensionalization of the compressible Navier-Stokes equations, see for instance \cite{PlotnikovSokolowski12}). Therefore, the small fluid velocity competes with the small length scale and local incompressibility can only be expected if $\sigma_\eps^2 \eps^{-\alpha} \to 0$, that is $\alpha < 3/2$.
We summarize this in the following conjecture.

\begin{conj} \label{conj}
Let $1 < \alpha < 3/2$, $\beta = 0$ and $\gamma > 3/2$. Let $(\ueps, \rhoeps)$ be a sequence of weak energy solutions of \eqref{goal}. Then the extensions by $0$, $(\tilde \rho_\eps, \tilde u_\eps)$ converge in a suitable way to $(\rho, u)$ that solve
\begin{align}
\begin{cases}
\partial_t \rho + \div \left( \rho u  \right)=  0 & \text{ in } \Omega \times (0,T).\\
u = \Rcal^{-1} (\rho f + g - \nabla \rho^{\gamma}) \text{ on } \{ \rho > 0 \} & \text{ in }\Omega \times (0,T), \\
u \cdot n = 0 &\text{ on } \partial \Omega \times (0,T),
\end{cases}
\end{align}
\end{conj} 
As a side remark, this rises the intriguing question about the homogenization limit for $\alpha \in [3/2,3]$. 

The conjecture seems a highly challenging problem. One of the main difficulties appears to be the lack of regularity of $\rho_\eps$ which makes it difficult to catch the ``almost incompressibility'' of $u_\eps$ around the particles.
Indeed, as we will see in the proof of our result, we will need our assumption on $\beta > 0$ precisely at the point when we need to prove smallness of terms involving $\div u_\eps$. 

One can view our result as a step towards  Conjecture \ref{conj}.
We show that for $\gamma \geq 2$, we only need to require $\beta > (3/2) (\alpha - 1) (\gamma + 2)$ (see Theorem \ref{theorem3} below). Thus, we can treat any $\beta > 0$ if $(\alpha - 1)$ is sufficiently small. This relates to the result by Masmoudi \cite{Masmoudi} for $\alpha=1$, $\beta = 0$.

\smallskip 

We emphasize though that the significance of our result goes beyond this relation to  Conjecture \ref{conj}. Indeed, in many applications the Mach number is very small motivating the study of low Mach limits. The low Mach number limit of compressible fluids has  been rigorously studied on a fixed smooth domain in the seminal papers by Lions, Masmoudi, Desjardins and Grenier in \cite{lowMach_98} and \cite{lowMach_99}.
However, to our knowledge, we study here for the first time the low Mach limit for a homogenization problem.

Moreover, our result covers all adiabatic exponents $\gamma > 3/2$ which is precisely the regime where existence of weak energy solutions to the compressible Navier-Stokes system \eqref{goal} is known. In the previously cited homogenization results 
the exponents are restricted to $\gamma > 2$ (\cite{FeireislNovotnyTakahashi10,LuPokorny20,DieFeiLu17}), $\gamma \geq 3$ (\cite{Masmoudi,FeireislLu15} and even $\gamma > 6$ (\cite{schwarzacher}).

\smallskip

Our proofs proceed in two steps, namely uniform a priori estimates and the passage to the limit in the weak formulation of the equation.

The basis of the second step lies in the classical method of oscillating test functions that goes back to Tartar \cite{tartar}, Cioranescu and Murat \cite{CioranescuMurat82a} and 
Allaire \cite{allaire2, allaire3}. Roughly speaking, the idea is  to use that the Stokes operator is self-adjoint in order to make rigorous the heuristics explained in Section \ref{sec:heuristics} on the level of an oscillating test function that can be constructed explicitly.
We will explain this in more detail in Section \ref{sec:testfunctions}.

For the first step, uniform a priori estimates for the fluid velocity follow from the energy inequality and a well-known quantitative Poincar\'e inequality in $\Omega_\eps$. The main difficulty lies in establishing uniform a priori estimates for the pressure. Here lies a technical novelty of our proof.
Following the classical approach due to Tartar one can define a suitable restriction operator $R_\eps \colon H^1_0(\Omega) \to H^1_0(\Omega_\eps)$.
By duality this allows to establish estimates on the pressure through estimates on the restriction operator.
We do not follow this approach here.
Instead, we characterize the pressure by using the Bogovski\u{\i} operator $\Bcal_\e$ on the porous domain that has been studied in \cite{schwarzacher} (see Theorem~\ref{thm_bog_eps} below). The advantage of this methodology is that one has estimates in the full range of exponents of Lebesgue and Sobolev spaces with constants depending on $\e$ in the expected way which seems not to be known for the restriction operator. These estimates are needed since we treat all adiabatic exponents $\gamma > 3/2$.
In the literature, it is sometimes argued that using the restriction operator yields a (natural) way to extend the pressure inside the particles. Indeed, in the case $\alpha = 1$, the choice of the extension seems to be important. However, we just work with the extension by zero. Indeed, since $\alpha > 1$, the volume fraction of the particles vanishes as $\eps \to 0$, and thus the precise way to extend the pressure does not matter. Moreover, in our view, there seems to be no physical meaning of the choice of the pressure inside the particles.

\subsection{Organization of the paper}


The rest of the paper is organized as follows.

In Section \ref{sec:main.results}, we give the precise statements of our main results.
We will first state our result for a simplified stationary version of the compressible Navier-Stokes equations,
namely
\begin{align}\label{system1}
\begin{cases}
- \Delta \ueps + \frac{1}{\eps^\beta}\nabla \rhoeps^{\gamma} = \rhoeps f + g &  \text{ in }\omegaeps ,  \\
\div(\rhoeps \ueps)  = 0  &  \text{ in }\omegaeps,  \\
\ueps = 0  &\text{ on } \partial \omegaeps.
\end{cases} 
\end{align}
We include this result mainly for the sake of the presentation of its proof,
which is much more clear-cut than for the full system, yet containing most of the main ingredients. The convergence for the steady model holds for all $\gamma \geq 1$. 
We also believe that in order to attack Conjecture \ref{conj}, it could be convenient to start from this stationary model.

In Section \ref{sec:Prep}, we collect some important tools that will be used in the proof. These include the Poincar\'e inequality in $\Omega_\eps$, the Bogovski\u{\i} operator as well as the existence and properties of the oscillating test functions  used by Allaire. As we will rely on some estimates of these test functions in Sobolev spaces
which have not been provided by Allaire, we will also recall the construction of these test functions.

Finally, Sections \ref{sec:Stokes} and \ref{sec:Navier.Stokes} are devoted to the proof of the main results for the steady problem and the full problem respectively.
In both sections, we proceed similarly, first proving uniform a priori estimates and then passing to the limit.

\section{Statement of the main results} \label{sec:main.results}

\subsection{Notation}

For $T>0$, we introduce the notation $L^q_T(L^a(\Omega)) :=L^q(0,T;L^a(\Omega))$ 
and  $L_T^q(\Omega):= L^q_T(L^q(\Omega))$.

Since our solutions $(\ueps, \rhoeps)$ are defined in different domains $\omegaeps$, it is convenient to consider them as functions on the fixed reference domain $\Omega$. Therefore, we introduce the following extension by $0$.
For any function $h\in L^1(\omegaeps)$, we define \begin{equation}
\widetilde{h} := \begin{cases}
h & \text{in } \omegaeps \\
0 & \text{in } \Omega\setminus \omegaeps.
\end{cases}
\end{equation}
%

We introduce the following notation for the average over $\Omega_\eps$ of a function $h \in L^1(\Omega_\eps)$:
\[
\dpr{h}_{\e}=\frac{1}{\abs{\Omega_\e}}\int_{\Omega_e}h(y)\, dy.
\]

It is convenient to assume that the integral of the pressure vanishes. Since the equations are invariant under adding a constant to the pressure, we can achieve this by defining
\begin{align} \label{def:p}
	p_\eps := \frac{1}{\eps^\beta} \left(\rho_\eps^\gamma - \lrangle{\rho_\eps^\gamma}_\eps \right).
\end{align}

\subsection{Main result for the steady system}

It seems that  system \eqref{system1} has not been studied a lot in the literature.
However, since we essentially regard it as a toy model for the full equations, we just take for granted the existence of the following weak energy solutions.
Let $\gamma \geq 1$, $g \in L^2(\Omega, \R^3)$ and $f \in L^a(\Omega, \R^3)$, where $\frac 1 a + \frac 1 2 + \frac 1 {2\gamma} = 1$. Let $m_0 >0$ be independent of $\varepsilon$.
Then, we call $(\rho_\eps, u_\eps) \in L^{2 \gamma}(\Omega_\eps) \times H^1_0(\Omega_\eps, \R^3)$ a weak energy solution if $\rho_\eps \geq 0$ satisfies
$\int_{\Omega_\eps} \rho_\eps = m_0$, $(\rho_\eps, u_\eps)$ is a distributional solution to \eqref{system1} and satisfies 
the energy inequality
\begin{equation} \label{energy_eq_1}
\int_{\omegaeps} |\nabla \ueps|^2\, dx \leq  \int_{\omegaeps} \rhoeps f \cdot \ueps\, dx + \int_{\omegaeps} g \cdot \ueps\, dx.
\end{equation}

We remark that the energy equality can be formally derived as an equality by testing the equation with the solution.
More precisely, first we are testing the weak momentum equation with $\ueps$ 
\begin{equation} \label{testing}
\int_{\omegaeps} |\nabla \ueps|^2\, dx - \frac{1}{M_{\varepsilon}}\int_{\omegaeps} \rhoeps^{\gamma} \div \ueps\, dx = \int_{\omegaeps} f  \rhoeps \ueps + g \ueps\,dx.
\end{equation}
Then the following (formal) computation shows that the second integral on the left handside vanishes. Indeed, for $\gamma >1$,
\begin{gather*}
\int_{\omegaeps} (\nabla \rhoeps^{\gamma}) \ueps\, dx  = \int_{\omegaeps} \gamma \rhoeps^{\gamma - 2} (\nabla \rhoeps) (\rhoeps \ueps)\, dx = \frac{\gamma}{\gamma - 1} \int_{\omegaeps} \nabla \rhoeps^{\gamma - 1} (\rhoeps \ueps)\, dx\\
 = - \frac{\gamma}{\gamma-1}  \int_{\omegaeps} \rhoeps^{\gamma - 1} \div (\rhoeps \ueps)\, dx = 0,
\end{gather*}
and similarly, for $\gamma = 1$ we have
\[
\int_{\omegaeps} (\nabla \rhoeps) \ueps\, dx  = \int_{\omegaeps} \log(\rhoeps) \div(\rhoeps \ueps) = 0.
\]
Note that $\rhoeps$ is not in any Sobolev space so that the terms in the computation are not well defined.

%

We can now state the homogenization result for \eqref{system1}.

\begin{theorem}\label{theorem1}
Let $\alpha \in (1,3)$, $\gamma \geq 1$ and $f,g$ as above. Let $(\ueps, \rhoeps)$ be a sequence of weak energy solutions to system \eqref{system1}. Assume
 \[
\beta >\frac{3}{2}(\gamma + 1)(\alpha - 1).\]
Then, with $p_\eps$ defined as in \eqref{def:p},
\begin{alignat*}{2}
\rhotilde &\to  \rho_0  \quad && \text{ strongly in } L^{2 \gamma}(\Omega),\\
\tilde p_\e &\weak p \quad&& \text{ weakly in }  L^{2}(\Omega),\\
\frac{\utilde}{\sigmaeps^2} &\weak u \quad && \text{ weakly in } L^{2}(\Omega, \R^3),
\end{alignat*}
where $\rho_0 = m_0/|\Omega|$ and $(u,p) \in L^2(\Omega, \R^3) \times  H^1(\Omega)$ is the
unique weak solution with $\int_\Omega p = 0$ to the Darcy's law
\begin{equation}\label{limit_system1}
\begin{cases} u = \Rcal^{-1}(\rho_0 f + g - \nabla p) \quad & \text{in } \Omega , \\
\div u =  0 \hspace{1 ex} \quad & \text{in } \Omega  ,\\
u\cdot n = 0 \hspace{1 ex} \quad & \text{on } \partial \Omega  ,
\end{cases}
\end{equation} 
where $\Rcal$ is the resistance matrix defined in \eqref{def:Rcal}.
\end{theorem}

\begin{rem}
Note that the solution $p$ to \eqref{limit_system1} satisfies the Poisson equation with Neumann boundary data in $\Omega$. Therefore, $p$ is indeed unique up to additive constants. 
\end{rem} 


\subsection{Main result for the evolutionary compressible Navier-Stokes equations}

%
Existence of weak energy solutions for the  barotropic compressible Navier-Stokes equations \eqref{goal} has been proved in the seminal works of Lions~\cite{lions} for $\gamma \geq \frac 9 5$  and Feireisl, Novotn\'y and Petzeltova~\cite{feireisl} for $\gamma > \frac 3 2$ (see also the well-written book by Novotn\'y and Stra\v{s}kraba \cite{novotny}). We recall here the relevant parts of the existence result for future reference.

Let $T>0$, $f,g \in L_T^{\infty}(\Omega, \R^3)$, and assume
\begin{gather}
\rho_{\varepsilon 0} \in L^{\gamma}(\omegaeps), \hspace{2 ex} \rho_{\varepsilon 0} \geq 0 \text{ a.e. in }\omegaeps, \\
q_{\varepsilon 0} \in L^{\frac{2 \gamma}{\gamma + 1}}(\omegaeps, \R^3), \hspace{2 ex} q_{\varepsilon 0} = 0 \text{ on } \{ \rho_{\varepsilon 0} = 0\}, \\
\rho_{\varepsilon 0} | u_{\varepsilon 0} |^2 \in L^1(\omegaeps), \text{ where }  u_{\varepsilon 0} = \frac{q_{\varepsilon 0}}{\rho_{\varepsilon 0}} \text{ on } \{ \rho_{\varepsilon 0} >0 \} \text{ and }u_{\varepsilon 0} =0
\text{ on } \{ \rho_{\varepsilon 0} =0 \},
\end{gather}
Then, there exists $(\rho_\eps,u_\eps) \in L^\infty_T(\Omega_\eps) \times L^2_T(H_0^1(\Omega_\eps, \R^3))$ with $\rho_\eps \geq 0$, that solves \eqref{goal} in distribution. Moreover, the solution satisfies conservation of mass and the energy inequality, i.e, for all $t \in (0,T)$, 
\begin{equation}
\label{mass_conservation}
	\inteps \rhoeps(t)\, dx = \inteps \rho_{\varepsilon 0}\, dx
\end{equation}
and 
\begin{gather}
\frac{\sigmaeps^2}{\gamma - 1} \frac{1}{\varepsilon^{\beta}}\int_{\omegaeps} \left( \rhoeps^{\gamma}(t)- \rho_{\varepsilon 0}^{\gamma}\right)\, dx  + \sigmaeps^2 \int_{\omegaeps} \frac{\rhoeps(t)|\ueps(t)|^2}{2}\, dx + \int_0^t \int_{\omegaeps} |\nabla \ueps|^2 + (\div \ueps)^2\, dx\, dt \nonumber \\ 
\leq \sigmaeps^2  \int_{\omegaeps} \frac{\rho_{\varepsilon 0}| u_{\varepsilon 0}|^2}{2}\, dx + \int_0^t \int_{\omegaeps} (\rhoeps f + g) \cdot \ueps\, dx\, dt. \label{energy_eq}
\end{gather}
Furthermore, the following extended continuity equations holds in the sense of distributions:
\begin{align} \label{eq:continuity.everywhere.assumption}
	\sigmaeps^2 \partial_t \tilde \rho_\eps + \div(\tilde \rho_\eps \tilde u_\eps) = 0 \quad \text{in } \R^3 \times (0,T).
\end{align}

We impose the following additional assumption regarding the initial data, which imply that the density is close to constant. Such assumptions are typical for the study of low Mach limits. 
\begin{align} 
\lrangle{\rho_{\eps0}}_\eps \to \rho_0 \quad \text{for some } \rho_0 > 0, \\
\label{assumption_rho}
\frac{1}{\eps^\beta} \int_{\omegaeps} \left|  \rho_{\varepsilon 0}^\frac{\gamma}{2}  - \langle \rho_{\varepsilon 0} \rangle_\e^\frac{\gamma}{2} \right|^2\,dx \leq C,
\end{align}
where $C$ is independent of $\eps$.


We are now ready to state our main result.
\begin{theorem}\label{theorem3}
Let $\alpha \in (1,3)$, $\gamma > \frac{3}{2}$ and $f,g,\rho_{\eps0}, q_{\eps0}$ as above and let  $(\rho_\eps,u_\eps) \in L^\infty_T(\Omega_\eps) \times L^2_T(H_0^1(\Omega_\eps, \R^3))$ 
be a weak solution to \eqref{goal} with the properties specified above.
Assume \[
\beta >\frac{3}{2}(\gamma + 2)(\alpha - 1).\]
If $\gamma<2$, assume additionally
\begin{align}
\label{eq:larger.beta}
\beta > \max \{ 2\gamma(3-\alpha), \; 3(\alpha - 1) + \gamma(\alpha+3)  \}.
\end{align}
Then there exists a decomposition of the pressure $p_\eps$, defined as in \eqref{def:p},
\[
\tilde p_\e=p_{\e,1}+p_{\e,2},
\]
such that 
\begin{alignat*}{2}
\rhotilde^\frac{\gamma}{2} &\to  \rho_0^\frac{\gamma}{2}  \quad && \text{ strongly in } L^{\infty}(0,T;L^{2}(\Omega)),
\\
\frac{\utilde}{\sigmaeps^2} &\weak u \quad && \text{ weakly in } L^{2}([0,T]\times\Omega, \R^3),
\\
p_{\e,1}&\to 0 \quad &&\text{ strongly in } H^{-1}(0,T;L^q(\Omega)), \\
p_{\e,2}&\weak p \quad &&\text{ weakly in }L^2([0,T]\times \Omega) 
\end{alignat*}
for some $q\in (1,2]$,
%
and $(u,p) \in L^2_T(\Omega) \times L^2_T(H^1(\Omega, \R^3))$ is the unique solution 
to the Darcy's law \eqref{darcy_law1} with $\int_\Omega p(t) = 0$ for all $t \in (0,T)$.
\end{theorem}
\begin{rem}
In the decomposition of the pressure, $p_{\eps,1}$ is the pressure corresponding to inertial forces, and $p_{\eps_2}$ is the pressure corresponding to viscous and external forces.	The decomposition  will be made precise in \eqref{eq:p1} and \eqref{eq:p2}. 
\end{rem}

\section{Preparations} \label{sec:Prep}

In this section, we collect some known results which have already been obtained and used in the study of homogenization for fluids. 
These include a strengthened Poincaré inequality for the family of perforated domains $\{ \omegaeps\}$, estimates on the Bogovski\u{ı} operator in these domains, as well as
the properties and estimates of the oscillating test functions introduced by Allaire.
We emphasize that it is crucial for the analysis of the homogenization limit to analyze the (optimal) dependence of the estimates on $\eps$.

%

Throughout this paper, we adopt the notation that $C$ denotes a constant which is independent of $\varepsilon$. The exact values of these constants might differ from line to line. 
From now on, to shorten the notation, we will not distinguish between scalar and vector fields in the notation of Lebesgue and Sobolev spaces, i.e. $\phi \in W^{k,r}(\Omega)$ can either be scalar or a vector field, which will be clear from the context.

\subsection{Poincaré's inequality  and \texorpdfstring{Bogovski\u{ı}}{Bogovski} operator in the perforated domain}

The following Poincar\'e inequality in $\Omega_\eps$ taken from \cite{allaire2} will be crucial for our analysis.

\begin{lem}[{\cite[Lemma 3.4.1]{allaire2}}] \label{poincare}
	There exists a constant $C$ independent of $\varepsilon$ such that for all $u \in H^1_0(\omegaeps)$ 
	\begin{equation} \label{poicare_ineq}
	\| u \|_{L^2(\omegaeps)} \leq C \sigmaeps \| \nabla u \|_{L^2(\omegaeps)},
	\end{equation}
	where $\sigmaeps$ is defined as in $(\ref{sigma})$.
\end{lem}

An important tool for the study of the Navier-Stokes equations will be the Bogovski\u{ı} operator. We define the space $L_0^p(\Omega)$ by \[
L_0^p(\Omega) = \left\{ f \in L^p(\Omega) : \int_{\Omega} f \, dx =0\right\}.
\]
%


The Bogovski\u{\i} operator we use has originally be introduced in~\cite{DieRuzSch10} and was addapted to homegeniously perforated domains $\omegaeps$ in~\cite{DieFeiLu17} It was than later refined in \cite{schwarzacher} (in order to be able to treat the  unsteady setup). We will built our argument on this refinement, namely on~\cite[Proposition~2.2]{schwarzacher}:

\begin{theorem}\label{thm_bog_eps}
Let $\omegaeps$ be defined as in \eqref{omegaeps} with $\alpha \geq 1$, then for any $1<a<\infty$, there exists a linear operator $\mathcal{B}_{\varepsilon} : L_0^a(\omegaeps) \to W_0^{1,a}(\omegaeps)$ such that 
\begin{equation}\label{propbogeps}
\div\mathcal{B}_{\varepsilon}(f) = f \text{ in } \omegaeps, \hspace{3 ex} \| \mathcal{B}_{\varepsilon}(f) \|_{W_0^{1,a}(\Omega_\eps)} \leq C \left(1+ \varepsilon^{ -\alpha  + 3 \frac{\alpha - 1}{a}} \right)  \| f \|_{L^a(\omegaeps)}
\end{equation}
for some constant $C$ independent of $\varepsilon$.
For any $q > 3/2$, the linear operator $\mathcal{B}_{\varepsilon}$ can be extended as a linear operator from $\{ \div g: g \in L^q(\omegaeps), \, g \cdot n = 0 \text{ on } \partial \omegaeps \}$ to $L^q(\omegaeps)$ satisfying \begin{equation}\label{negative bog}
\| \mathcal{B}_{\varepsilon}(\div g)\|_{L^q(\omegaeps)} \leq C \| g \|_{L^q(\omegaeps, \mathbb{R}^3)},
\end{equation}
for some constant $C$ independent of $\varepsilon$.
\end{theorem}

\begin{rem}
In particular, for $a=2$ and $\alpha \in [1,3]$, the estimate in \eqref{propbogeps} reads \begin{equation}
\label{eq:bog2}
\| \mathcal{B}_{\varepsilon}(f) \|_{H_0^{1}(\Omega)} \leq C \left(1+ \varepsilon^{\frac{\alpha-3}{2}} \right)  \| f \|_{L^2(\omegaeps)} \leq  \frac{C}{\sigmaeps} \| f \|_{L^2(\omegaeps)}.
\end{equation}
\end{rem}

\subsection{The framework of  oscillating test functions} \label{sec:testfunctions}

As mentioned in the introduction, our convergence proof relies on the use of the oscillating test functions that have been constructed for the perforated domain $\Omega_\eps$ by Allaire in \cite{allaire1,allaire2}. We summarize the properties of these test functions that have been proved by Allaire in the following proposition.

\begin{prop}[{\cite[Proposition 3.4.12]{allaire2}}] \label{hypotheses}
Let $1< \alpha < 3$. Then, for $1 \leq k \leq 3$, there exist functions $(\wk, \qk)$ 
 such that\\[0.2 cm]
$(H1) \hspace{0.3 cm} \wk \in H^1(\Omega), \hspace{0.3 cm} \qk \in L^2(\Omega)$\\[0.1 cm]
$(H2) \hspace{0.3 cm} \div \wk = 0$ in $\Omega \hspace{0.2 cm}$ and $\hspace{0.2 cm} \wk = 0$ on the holes $T_i^{\varepsilon}$,\\[0.1 cm]
$(H3) \hspace{0.3 cm}\wk \to e_k$ in $L^2(\Omega)$, \hspace{0.3 cm} $\|\wk\|_{L^\infty(\Omega)} +  \sigmaeps \|\nabla \wk \|_{L^2(\Omega)} + \sigmaeps \|\qk \|_{L^2(\Omega)} \leq C$, \\[0.1 cm]
$(H4) \hspace{0.3 cm}$for all  sequences $\nu_{\varepsilon} \in H^1(\Omega)$ with
$\nu_{\varepsilon} \weak \nu \text{ in } L^2(\Omega)$  for some $\nu \in L^2(\Omega)$, 
$
\| \nabla \nu_{\varepsilon} \|_{L^2(\Omega)} \leq C/\sigmaeps$\\ \phantom{s} \hspace{0.75 cm}  and $
\nu_{\varepsilon} = 0$  in $\Omega \setminus \Omega_\eps$
and for all $\phi \in \mathcal{D}(\Omega)$,  \[
\sigmaeps^2 \langle \nabla \qk - \Delta\wk, \phi \nu_{\varepsilon} \rangle_{H^{-1}, H^1_0(\Omega)} \to  \int_{\Omega} \phi \, \Rcal e_k \cdot \nu \, dx,
\]
\phantom{s} \hspace{0.75 cm} where $\Rcal$ denotes the resistance matrix defined in \eqref{def:Rcal}. 
%
\end{prop}

%


Fur our purposes, we also need estimates for $\wk$ and $\qk$ in more general Sobolev spaces (see Lemma~\ref{lem:qk} below). Therefore, we repeat the construction of these functions, which is described in \cite{allaire2}.

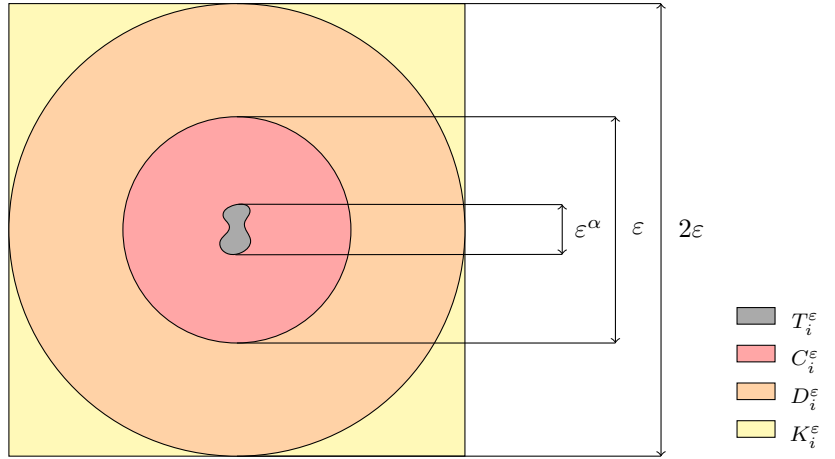
\begin{figure}[h] 
\hfill
\begin{tikzpicture} [scale = 0.1]
\filldraw[fill = yellow!35] (-27.8,-28.5) rectangle ++(60,60);
\filldraw[fill = orange!35] (2.2,1.5) circle (30);
\filldraw[fill = red!35] (2.2,1.5) circle (15);
\filldraw[fill={rgb:black,1;white,2}] (0,0) to [closed, curve through = {(1.2,2) (0.4,3.3) (1.9,4.8)  (3.9,4.2)(3.2,2.4) (4,0) (3,-1.4) (0.3,-1.2)}] (0,0);
 
\draw (2.4, 4.9) -- (45.0, 4.85);
 \draw (2.0, -1.8) -- (45.0, -1.75);
 \draw [<->] (45.0, 4.85) -- (45.0, -1.75);
 \node (C) at (48.5,1.7) {$\varepsilon^{\alpha}$};
 
 \draw (2.2, 31.5) -- (58.0, 31.5);
 \draw (2.2, -28.5) -- (58.0, -28.5);
  \draw [<->] (58.0, 31.5) -- (58.0, -28.5);
 \node (D) at (62.0,1.7) {$2 \varepsilon$};
 
 \draw [<->] (52.0, 16.5) -- (52.0, -13.5);
 \draw (2.2, 16.5) -- (52.0, 16.5);
 \draw (2.2, -13.5) -- (52.0, -13.5);
 \node (E ) at (55.0,1.7) {$ \varepsilon$};

 \filldraw[fill = yellow!35] (68,-26) rectangle ++(5,2.2);
 \filldraw[fill = orange!35] (68,-21) rectangle ++(5,2.2);
 \filldraw[fill = red!35] (68,-16) rectangle ++(5,2.2);
 \filldraw[fill={rgb:black,1;white,2}] (68,-11) rectangle ++(5,2.2);
 
 \node [font=\fontsize{8}{0}\selectfont](F) at (77.0,-26) {$K_i^{\varepsilon}$};
 \node [font=\fontsize{8}{0}\selectfont](G) at (77.0,-21) {{$D_i^{\varepsilon}$}};
 \node [font=\fontsize{8}{0}\selectfont] (H) at (77.0,-16) {{$C_i^{\varepsilon}$}};
 \node[font=\fontsize{8}{0}\selectfont] (I) at (77.0,-11) {{$T_i^{\varepsilon}$}};

\end{tikzpicture}
\caption{Decomposition of cell $\overline{P}_{\varepsilon}^{i}$}
\label{figure cell decomposition} 
\end{figure}

We decompose each cell, $\overline{P}^{\varepsilon}_{i}$, into the following four parts, illustrated in Figure \ref{figure cell decomposition}: \[
\overline{P}_{i}^{\varepsilon} = T_i^{\varepsilon} \cup \overline{C}_i^{\varepsilon} \cup \overline{D}_i^{\varepsilon} \cup \overline{K}_{i}^{\varepsilon}.
\]
Here $C_i^{\varepsilon}$ is the open ball of radius $\varepsilon/2$ centered in $P_i^{\varepsilon}$ and perforated by the hole $T_i^{\varepsilon}$; $D_i^{\varepsilon} = B_i^{\varepsilon} \setminus B_i^{\varepsilon/2}$ is the ball of radius $\varepsilon$ centered in $P_i^{\varepsilon}$ and perforated with a ball of radius $\varepsilon/2$ with the same center, and $K_i^{\varepsilon}$ are the remaining corners of the cell.

We define functions $(\wk, \qk)_{1 \leq k \leq 3} \in H^1(P_i^{\varepsilon}) \times L^2(P_i^{\varepsilon})$ with $\int_{P_i^{\varepsilon}} \qk =0$ by \[
\wk = e_k, \hspace{2 ex} \qk =0 \hspace{2 ex} \text{ in } P_i^{\varepsilon} \cap  \Omega
\]
for each cube $P_i^{\varepsilon}$  such that  $P_i^{\varepsilon} \cap \partial \Omega \neq \emptyset$ and by \[
\begin{cases}
\wk = e_k\\
\qk = 0
\end{cases}
\text{ in } K_i^{\varepsilon},
\hspace{1 cm}
\begin{cases}
\nabla \qk - \Delta \wk &=0\\
\div(\wk) &=0
\end{cases}
\text{ in } D_i^{\varepsilon},
\]\[
\begin{cases}
\wk &= w_k(\frac{x}{\varepsilon^{\alpha}})\\
\qk &= \frac{1}{\varepsilon^{\alpha}}q_k(\frac{x}{\varepsilon^{\alpha}})
\end{cases}
\text{ in } C_i^{\varepsilon},
\hspace{1 cm}
\begin{cases}
\wk &=0\\
\qk &= 0 
\end{cases}
\text{ in } T_i^{\varepsilon}
\]
for each cube entirely included in $\Omega$. Here the functions $w_k$ and $q_k$ are the solutions of the Stokes problem \eqref{stokessingle}, and the Stokes equations in $D^\eps_i$ are complemented with boundary conditions such that $w_k^\eps \in H^1(P_i^\eps)$.

\begin{lem} \label{lem:qk}
Let $\wk,\qk$ be defined as above.
Then, for all $p > \frac 3 2 $, there exists $C > 0$ independent of $\eps$ such that
\begin{align} 
	\label{eq:qk.C_i} \| \nabla w^\eps_k \|_{L^{p}(\Omega)}  +  \| \qk \|_{L^{p}(\Omega)}   &\leq C \eps^{-\alpha + 3\frac{\alpha - 1}{p}},  \\
	\label{eq:nabla.qk.C_i} \| \nabla \qk \|_{L^{p}\left(\cup_i C_i^{\varepsilon} \right)}  &\leq C  \eps^{-2\alpha + 3\frac{ \alpha - 1}{p}}.
\end{align}
Moreover, with $B^{\lambda}_i := B_\lambda(x_i^\eps)$,
\begin{align} 
\label{eq:qk.annuli}
	\| \nabla w_k^\eps \|_{L^2(\cup_i B_i^{\varepsilon}\setminus B_i^{\varepsilon /4})} +  \| \qk \|_{L^2(\cup_i B_i^{\varepsilon}\setminus B_i^{\varepsilon /4})} &\leq   C \varepsilon^{ \alpha - 2}.
\end{align}
\end{lem}

\begin{proof}
We first compute the norms of $q_k$ and $w^\eps_k$ in $\cup_i C_i^{\varepsilon}$.
To show \eqref{eq:qk.C_i} and \eqref{eq:nabla.qk.C_i}, we observe that by definition of $q_k^\eps$ in $C^\eps_i$ and scaling considerations we have
\begin{align}
	 \| \nabla w_k^\eps \|_{L^{p}\left(\cup_i C_i^{\varepsilon} \right)} + \| \qk \|_{L^{p}\left(\cup_i C_i^{\varepsilon} \right)} & \leq C \eps^{-\frac 3 p} \eps^{-\alpha + \frac{3 \alpha}{p}} \left(\|\nabla w_k\|_{L^p(\R^3\setminus T)} +\|q_k\|_{L^p(\R^3\setminus T)}\right), \\
	 \| \nabla \qk \|_{L^{p}\left(\cup_i C_i^{\varepsilon} \right)} & \leq C \eps^{-\frac 3 p} \eps^{-2\alpha + \frac{3 \alpha}{p}} \|\nabla q_k\|_{L^p(\R^3\setminus T)},
\end{align}
where $q_k$ is the pressure in the single particle problem \eqref{stokessingle} and the factor $\eps^{-\frac 3 p}$ arises because the number of particles is of order $\eps^{-3}$. By standard regularity theory for the Stokes equations 
(see e.g. \cite[Chapter V]{galdi}), $w_k$ and $q_k$ are smooth and satisfy for all $l \in \N$
\begin{align} 	\label{eq:decay.w_k}
\begin{aligned} 
	|w_k - e_k| &\leq C \frac{1}{|x|}, \\
 |\nabla^{l+1} w_k(x)| + |\nabla^l q_k(x)| &\leq C_k \frac{1}{|x|^{l+2}} \quad \text{in } \R^3 \setminus T. 
\end{aligned} 
\end{align}
In particular, $\nabla w_k \in L^p(\R^3 \setminus T)$ and $q_k \in W^{1,p}(\R^3 \setminus T)$ for all $p > \frac 3 2$.
This shows \eqref{eq:qk.C_i} and \eqref{eq:nabla.qk.C_i} once we established \eqref{eq:qk.annuli}.

\medskip
In order to prove \eqref{eq:qk.annuli}, we first estimate the norms in $D_i^\eps$. Again, by scaling considerations,
\begin{align}
	\| \nabla w_k^\eps \|_{L^{p}\left(\cup_i D_i^{\varepsilon} \right)} + \| \qk \|_{L^{p}\left(\cup_i D_i^{\varepsilon} \right)}   \leq C \eps^{-1} \left(\|\nabla v^\eps_k\|_{L^p(B_1(0) \setminus B_{1/2}(0))} +  \|p^\eps_k\|_{L^p(B_1(0) \setminus B_{1/2}(0))} \right),
\end{align}
where $v_k^\eps, p_k^\eps$ is the solution to homogeneous Stokes equations in $B_1(0) \setminus B_{1/2}(0)$ with boundary data
\begin{align}
	v_k^\eps = \begin{cases} 0 &\quad \text{on } \partial B_1(0), \\
	w_k\left( \frac {\eps \cdot}{2 \eps^\alpha}\right) - e_k &\quad \text{on } \partial B_{1/2}(0).
	\end{cases}
\end{align}
By standard regularity theory of the Stokes equations,
\begin{align} 
	\|\nabla v^\eps_k\|_{L^p(B_1(0) \setminus B_{1/2}(0))} +  \|p^\eps_k\|_{L^p(B_1(0) \setminus B_{1/2}(0))} 
	&\leq C \left\|w_k\left( \frac {\eps \cdot}{2 \eps^\alpha}\right) - e_k \right\|_{W^{1-\frac 1 p,p}(\partial B_{1/2}(0))}  \\
	& \leq C \|\nabla \phi\|_{L^p(B_1(0) \setminus B_{1/2}(0))},
\end{align}
for any $\phi \in W_{1,p}(B_1(0) \setminus B_{1/2}(0))$ satisfying the same boundary conditions as $v_k^\eps$. 
Such a function is given by $\phi = \eta \left(w_k\left( \frac {\eps \cdot}{2 \eps^\alpha}\right) - e_k\right)$ where $\eta$ is a suitable cutoff function with $\eta = 1$ on $\partial B_{1/2}(0)$ and $\eta = 0$ on $B_1(0)$. 
Then $\|\nabla \phi\|_\infty \leq \eps^{\alpha - 1}$ due to \eqref{eq:decay.w_k}. Thus,
\begin{align}
	\| \nabla w_k^\eps \|_{L^{p}\left(\cup_i D_i^{\varepsilon} \right)} + \| \qk \|_{L^{p}\left(\cup_i D_i^{\varepsilon} \right)} \leq C \eps^{\alpha -2}.
\end{align}
It remains to estimate $\nabla w_k^\eps$ and $q_k^\eps$ on $B^{\eps/2}_i \setminus B^{\eps/4}_i$. A straightforward computation using the decay estimates \eqref{eq:decay.w_k} yields 
\begin{align}
\| \nabla w_k^\eps \|_{L^{p}\left(\cup_i B^{\eps/2}_i \setminus B^{\eps/4}_i \right)} + \| \qk \|_{L^{p}\left(\cup_i B^{\eps/2}_i \setminus B^{\eps/4}_i \right)} \leq C \eps^{\alpha -2}.
\end{align}
This establishes \eqref{eq:qk.annuli}.
\end{proof}

\subsection{Auxiliary estimates}
We collect here some equivalence of quantities in the sense that each member can be estimated by the other with a uniform constant. For that we introduce the symbol ''$\sim$''.
We will use the following inequality, that can be found for instance in~\cite[Lemma 2.2]{GiaSch19}:
Let $a\in (\frac{1}{2},\infty)$.
for $f\in L^{2a}(\Omega_\e;[0,\infty))$, we find
\begin{align}
\label{eq:gamma1}
\int_{\Omega_\eps}\abs{f^a-\dpr{f}_\eps^a}^2\, dx\sim  \int_{\Omega_\eps}\abs{f^a-\dpr{f^a}_\eps}^2\, dx.
\end{align}
Observe that we find for $a,b\in [0,\infty)$ and $\gamma>1$ that
\[
(a+b)^{\gamma-2}\sim \int_0^1(b+\theta(a-b))^{\gamma-2}\, d\theta
\]
this implies for one that
\[
(a^{\gamma-1}-b^{\gamma-1})(a-b)\sim (a+b)^{\gamma-2}(a-b)^2\sim (a^\frac{\gamma}{2}-b^\frac{\gamma}{2})^2,
\]
 but also 
\begin{align}
\label{eq:gamma2}
\begin{aligned}
\frac{a^\gamma}{\gamma}-\frac{b^\gamma}{\gamma}-b^{\gamma-1}(a-b)&=\int_0^1\int_0^1(\gamma-1)(b+\theta\tau(a-b))^{\gamma-2}\, d\theta \tau(a-b)^2\, d\tau \\ &\sim (a+b)^{\gamma-2}(a-b)^2 
\sim (a^\frac{\gamma}{2}-b^\frac{\gamma}{2})^2.
\end{aligned}
\end{align}

\section{Homogenization of the steady compressible Stokes system including a Low-Mach-Number-Limit} \label{sec:Stokes}

\subsection{A priori estimates}

\begin{lem} \label{aprioriMach1}
Let the assumptions of Theorem \ref{theorem1} hold.
Then, 
 \begin{equation}\label{a-priori u 1}
\left\|\frac{\utilde}{\sigmaeps^2}\right\|_{L^2(\Omega)} + \left\|  \frac{\nabla \utilde}{\sigmaeps} \right\|_{L^2(\Omega)} \leq C,
\end{equation}
\begin{equation}\label{a-priori rho 1}
\|\rhoeps - \langle \rhoeps \rangle_\e\|_{L^{2 \gamma}(\omegaeps)} \leq C \eps^{\frac{\beta}{\gamma}},
\end{equation}
 and
\begin{equation}\label{a-priori p 1}
\| p_\e \|_{L^2(\Omega)} \leq C.
\end{equation}
\end{lem}

\begin{proof}


We start from the energy inequality \eqref{energy_eq_1}.
Using Hölder's inequality and Lemma \ref{poincare}, we estimate the right-hand side of \eqref{energy_eq_1} in the following way:

\begin{align*} \label{integral}
\int_{\omegaeps}\rhoeps f  \cdot \ueps\, dx + \int_{\omegaeps} g \cdot \ueps\, dx & \leq \|\ueps\|_{L^2(\omegaeps)}\left(\|f\|_{L^{a}(\omegaeps)} \| \rhoeps \|_{L^{2\gamma}(\omegaeps)} + \|g\|_{L^2(\omegaeps)} \right)\\
& \leq C \sigmaeps \| \nabla \ueps \|_{L^2(\omegaeps)} \left( \|f\|_{L^{a}(\omegaeps)} \| \rhoeps \|_{L^{2\gamma}(\omegaeps)} + \|g\|_{L^2(\omegaeps)} \right)\\
\end{align*}
so that from  $(\ref{energy_eq_1})$ we deduce 

\begin{equation}\label{uestimate}
\| \nabla \ueps\|_{L^2(\omegaeps)} \leq C \sigmaeps \left( \|f\|_{L^{a}(\omegaeps)} \| \rhoeps \|_{L^{2\gamma}(\omegaeps)} + \|g\|_{L^2(\omegaeps)} \right).
\end{equation}

We see that in order to get an estimate for $\ueps$, we need an estimate for $\rhoeps$ first. This will result from an estimate for the pressure $\peps$.
To this end, we use Theorem \ref{thm_bog_eps} to test the first equation in \eqref{system1} with $\Bcal_\eps(p_\eps) \in H^1_0(\Omega_\eps)$:
\begin{align}
\label{eq:defp2}
\| p_\eps\|_{L^2(\Omega_\eps)}^2  = \int_{\omegaeps} \peps \div \Bcal_\eps(p_\eps) \, dx = \int_{\omegaeps}\nabla \ueps : \nabla \Bcal_\eps(p_\eps)-(\rhoeps f+g)\cdot \Bcal_\eps(p_\eps) \, dx.
\end{align}
By H\"older's inequality  and 
\eqref{uestimate} we deduce with $a = \frac{2 \gamma}{\gamma - 1}$
\begin{align*}
\| p_\eps\|_{L^2(\Omega_\eps)}^2&\leq  \norm{\nabla \ueps}_{L^2(\omegaeps)}\norm{\nabla \Bcal_\eps(p_\eps)}_{L^2(\omegaeps)}+(\norm{f}_{L^a(\omegaeps)} \norm{\rhoeps}_{L^{2\gamma}(\omegaeps)}+\norm{g}_{L^2(\omegaeps)})\norm{\Bcal_\eps(p_\eps)}_{L^2(\omegaeps)}
\\
&\leq  C\left( \|f\|_{L^a (\omegaeps)} \| \rhoeps \|_{L^{2\gamma}(\omegaeps)} + \|g\|_{L^2(\omegaeps)} \right) (\sigmaeps \norm{\nabla \Bcal_\eps(p_\eps)}_{L^2(\omegaeps)} + \norm{\Bcal_\eps(p_\eps)}_{L^2(\omegaeps)}) \\
&\leq  C\left( \|f\|_{L^a (\omegaeps)} \| \rhoeps \|_{L^{2\gamma}(\omegaeps)} + \|g\|_{L^2(\omegaeps)} \right) \| p_\eps\|_{L^2(\Omega_\eps)},
\end{align*}
where we used Lemma \ref{poincare} and Theorem \ref{thm_bog_eps} in the last estimate.

Thus,
\begin{align}
\label{eq:pressure1}
\norm{p_\e}_{L^2(\Omega_\eps)}\leq C \left( \|f\|_{L^a (\omegaeps)} \| \rhoeps \|_{L^2(\omegaeps)} + \|g\|_{L^2(\omegaeps)} \right), 
\end{align}
with a constant independent of $\e$.
Now observe that \[
\langle \rhoeps \rangle_\e = \frac{1}{|\omegaeps|}\int_{\omegaeps} \rhoeps = \frac{m_0}{|\omegaeps|}.
\]
By \eqref{eq:gamma1} we find that
\begin{gather}
\frac{1}{\eps^\beta} \|\rhoeps^{\gamma} -\langle \rhoeps \rangle_\e^{\gamma}\|_{L^2(\omegaeps)} \leq C \frac{1}{\eps^\beta} \|\rhoeps^{\gamma} -\langle \rhoeps^{\gamma} \rangle_\e\|_{L^2(\omegaeps)} = C \|\peps\|_{L^2(\omegaeps)}.
\label{pres}
\end{gather}

Then, from \eqref{eq:pressure1} and \eqref{pres}, we find
\begin{align} \label{esimate_rho_gamma}
\frac{1}{\eps^\beta} \|\rhoeps^{\gamma} -\langle \rhoeps \rangle_\eps^{\gamma}\|_{L^2(\omegaeps)} &\leq C \|p_\e \|_{L^2(\Omega)}
 \leq C \left(\|f\|_{L^a (\omegaeps)} \|\rhoeps \|_{L^{2 \gamma}(\omegaeps)}  + \|g\|_{L^2(\omegaeps)} \right) \nonumber \\
& \leq C \left(\|f\|_{L^a (\omegaeps)} \|\rhoeps^\gamma-\dpr{\rhoeps}_\e ^\gamma\|_{L^{2}(\omegaeps)}^\frac{1}{\gamma}  +\frac{m_0\|f\|_{L^a (\omegaeps)}}{\abs{\omegaeps}^{1 - 1/(2 \gamma)}}+ \|g\|_{L^2(\omegaeps)} \right)
\\
& \leq  \frac{1}{2 \eps^\beta} \| \rhoeps^{\gamma} - \langle \rhoeps \rangle_\eps^{\gamma}\|_{L^2(\omegaeps)} +  C,
\end{align}
where we applied Young's inequality in the last inequality and used that we can assume  $\eps \leq 1$. 
Now using that $ | \rhoeps - \langle \rhoeps \rangle_\eps |^{ \gamma} \leq  | \rhoeps^{\gamma} - \langle \rhoeps \rangle_\eps^{\gamma} |$ (which follows for instance from the triangle inequality of the metric $d(a,b) = |a-b|^{1/\gamma}$), we conclude
\begin{gather*}
\frac{1}{\eps^\beta}\|\rhoeps - \langle \rhoeps \rangle_\eps \|_{L^{2 \gamma} (\omegaeps)}^{\gamma}  \leq  \frac{C}{\eps^\beta}\|\rhoeps^{\gamma} - \langle \rhoeps \rangle_\eps^{\gamma}\|_{L^2(\omegaeps)}\leq C
\end{gather*}
and hence we find \eqref{a-priori rho 1}.
Moreover,
 \[
\| \rhoeps\|_{L^{2 \gamma}(\omegaeps)} \leq \| \rhoeps - \langle \rhoeps \rangle_\e\|_{L^{2 \gamma}(\omegaeps)} + C \langle\rhoeps \rangle_\e \leq C.
\]
Finally, \eqref{uestimate} and Lemma~\ref{poicare_ineq}
imply \eqref{a-priori u 1} and \eqref{eq:pressure1} implies \eqref{a-priori p 1}.
\end{proof}

\subsection{Convergence Proof}

\begin{proof}[Proof of Theorem \ref{theorem1}]

\emph{Step 1:} Extracting (weakly) convergent subsequences.

We use the a priori estimates to find  $u \in L^2(\Omega), p \in L^2_0(\Omega)$ and  subsequences (not relabeled) such that
\begin{alignat*}{2}
\rhotilde &\to  \rho_0   \quad && \text{ strongly in } L^{2 \gamma}(\Omega),\\
\tilde p_\eps &\weak p \quad &&  \text{ weakly in }  L^{2}(\Omega),\\
\frac{\utilde}{\sigmaeps^2} & \weak u  \quad && \text{ weakly in } L^{2}(\Omega),
\end{alignat*}
where $\rho_0 = m_0/|\Omega|$ is constant.

To deduce the strong convergence of $\tilde \rho_\eps$, we estimate
\begin{align*}
\|\tilde \rhoeps - \rho_0 \|_{L^{2 \gamma} (\Omega)} & \leq \|\rho_0\|_{L^{2 \gamma}(\Omega \setminus \Omega_\eps)} +  \|\rhoeps - \langle \rhoeps \rangle_\eps\|_{L^{2 \gamma} (\omegaeps)} + \|\langle \rhoeps \rangle_\eps - \rho_0 \|_{L^{2 \gamma} (\omegaeps)} \\[0.1 cm]
& \leq C \eps^{\frac{\beta}{\gamma}} + |\omegaeps|^{\frac{1}{2 \gamma}} \left( \frac{1}{|\omegaeps|} - \frac{1}{|\Omega|} \right) m_0  + |\Omega \setminus \omegaeps|^{\frac 1 {2\gamma}} \rho_0 \to 0
\end{align*}
since $| \omegaeps| \to |\Omega|$.

It remains to prove that $(u,p)$ is a solution to Darcy's law \eqref{limit_system1}. Indeed, since this system has a unique weak solution, the convergence of the full sequence then follows.

\medskip

\emph{Step 2:} Passage to the limit in the continuity equation.

We pass to the limit in the the second equation of system \eqref{system1}. To this end, we prove that
\begin{align} \label{eq:divergencefree.R^3}
	\div (\tilde \rho_\eps \tilde u_\eps) = 0 \quad \text{in } \R^3.
\end{align}
We first show how to pass to the limit in the continuity equation assuming \eqref{eq:divergencefree.R^3} holds.
For $\phi \in \Dcal(\R^3)$, we have
\begin{align}
	0 = \int_{\R^3} \tilde \rho_\eps  \frac{\tilde u_\eps}{\sigma_\eps^2} \cdot \nabla \phi \to \int_\Omega u \cdot \nabla \phi,
\end{align}
which is the weak formulation of
\begin{align} \label{eq:u.divergencefree.1}
\begin{cases}
	\div u = 0 \quad & \text{in } \Omega, \\
	u \cdot n = 0 \quad & \text{on } \partial \Omega.
\end{cases}
\end{align}
 It remains to prove \eqref{eq:divergencefree.R^3}. Here, we follow the argument of \cite[Lemma 2.3]{Masmoudi} which we include for completeness.
Fix $\eps > 0$, and for $d \ll \eps$ let $\psi_{d} \in \Dcal(\R^3)$ be a family  of cutoff functions such that $0 \leq \psi_d \leq 1$,
$\psi_{d} = 0$ in $\R^3 \setminus \Omega_\eps$, $\psi_{d}(x) = 1$
for all $x \in \Omega$ with $\dist(x,\partial \Omega_\eps) \geq d$  and
$\|\nabla \psi_{d}\|_\infty \leq C d^{-1}$. In particular, $|\nabla \psi_d(x)| \leq C \dist^{-1}(x,\partial \Omega_\eps)|$. Hence, by Hardy's inequality
\begin{align}
	\|u_\eps \nabla \psi_d\|_{L^2(\Omega_\eps)} \leq C \|\nabla u_\eps\|_{L^2(\Omega)}.
\end{align}
Since,  $\psi_{d} \to \chi_{\Omega_\eps}$ as $d \to 0$ strongly in $L^p(\R^3)$ for every $p < \infty$, we deduce for all $\phi \in \Dcal(\R^3)$
\begin{align}
	0 = \int_\Omega \tilde \rho_\eps \tilde u_\eps \cdot \nabla (\psi_d \phi)\, dx
	= \int_\Omega \tilde \rho_\eps \tilde u_\eps \cdot \nabla \phi  \psi_d\, dx
	+ \int_\Omega \tilde \rho_\eps \tilde u_\eps \cdot \nabla \psi_d \phi\, dx
	\to \int_\Omega \tilde \rho_\eps \tilde u_\eps \cdot \nabla \phi\, dx
\end{align}
as $d \to 0$, where we used that
\begin{align}
	\left|\int_\Omega \tilde \rho_\eps \tilde u_\eps \cdot \nabla \psi_d \phi \right| \leq C \| \rho_\eps \chi_{\dist(x,\partial \Omega_\eps) \leq d}\|_{L^2(\Omega_\eps)} \|\nabla u_\eps \|_{L^2(\Omega_\eps)} \to 0.
\end{align}


\medskip

\emph{Setp 3:} Passage to the limit in the momentum equation.

Let $(\wk, \qk)$ be as in Proposition \ref{hypotheses}. Let $\phi \in \mathcal{D}(\Omega)$. Using  $\wk \phi$ as a test function in the first equation of \eqref{system1} yields (observe that$ \int_{\omegaeps}\div(\wk \phi)\, dx=0$)
\begin{equation}
\int_{\Omega} \nabla \utilde  :\nabla (\wk \phi)\, dx =  \int_{\Omega}  \tilde p_\e \div(\wk \phi)\, dx +  \int_{\Omega} (\rhotilde f + g)\cdot(\wk \phi)\, dx.
\end{equation}
We can integrate over the whole domain $\Omega$ since $\utilde$, $\wk$  vanish on $\Omega \setminus \omegaeps$. Since we want to make use of $(H4^*)$ we rewrite the left hand side in the following way:
\begin{align*}
\int_{\Omega} \nabla \utilde :\nabla (\wk \phi)\, dx  &=  \int_{\Omega} \nabla \utilde : \left(  \wk \otimes\nabla\phi + \nabla \wk \phi\right)\, dx \nonumber \\
 &= \int_{\Omega} \nabla \utilde : (\wk \otimes \nabla \phi)\, dx - \int_{\Omega} \nabla \wk : ( \utilde \otimes \nabla \phi)\, dx +  \int_{\Omega} \nabla \wk: \nabla\left(\utilde \phi\right)\, dx
\end{align*}
and add the term $ - \int_{\Omega}  \qk \div \left(\utilde \phi  \right)\, dx$ to both sides, which leads to 
\begin{gather}
\underbrace{ \int_{\Omega} \nabla \wk :  \nabla(\utilde \phi)- \qk \div \left(\utilde \phi \right)\, dx}_{I_1} + \underbrace{\int_{\Omega} \nabla \utilde:( \wk \otimes\nabla\phi ) - \nabla \wk : (\utilde \otimes \nabla\phi)\, dx}_{I_2}\nonumber \\
= \underbrace{\int_{\Omega} \tilde p_\eps \wk \cdot \nabla \phi +  (\rhotilde f + g) \cdot(\wk \phi)\, dx}_{I_3} -\underbrace{ \int_{\Omega} \qk \div \left(\utilde \phi \right)\, dx}_{I_4},
\end{gather}
where we used $\div(\wk) =0$ in the first integral on the right-hand side. 
By hypothesis $(H4)$ we then have
\begin{equation}
I_1 = \sigma_{\varepsilon}^2 \int_{\Omega} \nabla \wk \nabla \left(\frac{\utilde}{\sigma_{\varepsilon}^2} \phi \right)\, dx- \sigma_{\varepsilon}^2 \int_{\Omega} \qk \div \left(\frac{\utilde}{\sigma_{\varepsilon}^2} \phi \right)\, dx \to \int_{\Omega} \phi \Rcal e_k \cdot u\, dx
\end{equation}
as $\varepsilon \to 0$. Further by our a-priori estimates \eqref{a-priori u 1} and hypothesis $(H3)$,
\begin{equation}
|I_2| = \sigma_{\varepsilon} \left|\int_{\Omega} \left(\frac{\nabla  \utilde}{\sigma_{\varepsilon}} :(\wk \otimes \nabla \phi) - \sigma_{\varepsilon}\nabla \wk: \left( \frac{\utilde}{\sigmaeps^2}\otimes \nabla \phi\right) \right)\, dx \right|\leq C \sigmaeps\to 0 
\end{equation}
as $\varepsilon \to 0$. Since $\wk$ converges strongly to $e_k$ in $L^2(\Omega)$, we obtain 
\begin{equation}
I_3 \to  \int_{\Omega} \left(\rho_0 f +  g\right)\cdot e_k \phi + p e_k \cdot \nabla \phi\, dx.
\end{equation}
We postpone the convergence  $I_4 \to 0$ as $\varepsilon \to 0$ to the end of the proof.

Accepting this convergence for the moment, taking the limit in the momentum equation yields
\begin{equation}
 \int_{\Omega} \phi \Rcal e_k \cdot u \, dx =  \int_{\Omega} p \partial_k \phi + (\rho_0f + g)\cdot e_k \phi\, dx
\end{equation}
for all $\phi \in \mathcal{D}(\Omega)$. Since this holds true for all $1 \leq k \leq 3$ we infer that $\nabla p\in L^2(\Omega)$ and  
\begin{equation}\label{proof limit momentum equation 1}
\Rcal u = \rho_0 f + g - \nabla p.
\end{equation}
Combining \eqref{eq:u.divergencefree.1} with \eqref{proof limit momentum equation 1} yields \eqref{limit_system1}.

\medskip

\emph{Step 4:} Estimate for $I_4$.

To conclude the proof, we must show $I_4 \to 0$ as $\eps \to 0$.

First of all, observe that \begin{align*}
 \int_{\Omega} \qk \div(\utilde \phi)\, dx & =  \int_{\omegaeps} \qk \div(\ueps \phi)\, dx.
\end{align*}
Recall from Section \ref{sec:testfunctions} that in each cell $P_i^{\varepsilon}$ the function $\qk$ is defined piecewise on the sets $T_i^{\varepsilon}$, $C_i^{\varepsilon}$, $D_i^{\varepsilon}$ and $K_i^{\varepsilon}$, where we use the same notation as in Section \ref{sec:testfunctions}. We would like to integrate by parts in $I_4$, but by its piecewise construction, $\nabla \qk$ is not well defined as a function in some $L^a(\Omega)$. Nevertheless, it is welldefined inside each $C_i^{\varepsilon} = B_i^{\varepsilon/2} \setminus T_i^{\varepsilon}$, where here and in the following we write $B_i^{r_{\varepsilon}} = B_{r_{\varepsilon}}(x_i^{\varepsilon})$ for the ball of radius $r_{\varepsilon}$ centered at $x_i^{\varepsilon}$.

Therefore, the idea is to split the integral over $\omegaeps$ in the following manner:

Let $\psi \in C_c^{\infty}\left(\bigcup_i B_i^{\varepsilon/2}\right)$ be a cutoff function with \[
\psi = 1 \text{ in }  \bigcup_i B_i^{\varepsilon/4}, \hspace{2 ex} | \nabla \psi| \leq \frac{C}{\varepsilon},\]
and extend it by $0$ to $\Omega$.
Then, \[
\qk = \psi  \qk + (1-\psi) \qk,
\]
and we can write \begin{align}
\langle \rhoeps \rangle_\eps \cdot I_4 &= \langle \rhoeps \rangle_\eps \int_{\omegaeps} \qk \div(\ueps \phi)\, dx\nonumber \\
&=\langle \rhoeps \rangle_\eps \int_{\omegaeps} \qk \psi  \div(\ueps \phi) \, dx+ \langle \rhoeps \rangle_\eps \int_{\omegaeps} \qk (1 - \psi) \phi \div(\ueps) \, dx
 \nonumber \\
 & \hspace{0.5 cm}+ \sigmaeps^2 \langle \rhoeps \rangle_\eps \int_{\omegaeps} \qk (1 - \psi) \frac{\ueps}{\sigmaeps^2} \cdot\nabla \phi \, dx\nonumber \\
&=: I^1 + I^2 + I^3. \label{i1}
\end{align}
By hypothesis $(H3)$ and \eqref{a-priori u 1}, \[
|I^3| \leq  C \sigmaeps\to 0.
\]
In order to prove that $|I^1| \to 0$, we use that $\supp(\psi) \subset \cup_i B_i^{\varepsilon/2}$ and  $\div(\rhoeps \ueps) = 0$ (see \eqref{eq:divergencefree.R^3}). Thus,
\begin{align*} 
I^1 
& =   \int_{\omegaeps} \nabla (\qk \psi \phi) \rhoeps \ueps\, dx - \int_{\omegaeps} \nabla (\qk \psi \phi) \langle \rhoeps \rangle_\eps \ueps\, dx  +  \langle \rhoeps \rangle_\eps \int_{\omegaeps} \qk \psi \, \ueps \cdot \nabla \phi\, dx\\
&= \sigmaeps^2  \int_{\omegaeps} \nabla (\qk \psi \phi) \left(\rhoeps -\langle \rhoeps \rangle_\eps \right) \frac{\ueps}{\sigmaeps^2}\, dx + O(\sigmaeps),
\end{align*}
where for the last term we used hypotheses $(H3)$ and \eqref{a-priori u 1} in order to see that it converges to zero with order $\sigmaeps$.
Since $\supp (\psi) \subset \cup_i B_i^{\varepsilon/2}$ we find, recalling $C_i^\eps = B_i^{\varepsilon/2}\setminus T_i^{\varepsilon}$
 \begin{align*}
|I^1| & \leq \sigmaeps^2 \|\nabla (\qk \psi \phi)\|_{L^{\frac{2 \gamma}{\gamma - 1}}\left(\cup_i C_i^\eps \right)}  \| \rhoeps -\langle \rhoeps \rangle_\eps \|_{L^{2 \gamma}(\omegaeps)}  \left\| \frac{\ueps}{\sigmaeps^2} \right\|_{L^2(\omegaeps)} + O(\sigmaeps) \\
& \leq  C \sigmaeps^2 \eps^{\frac{\beta}{\gamma}} \|\nabla (\qk \psi \phi)\|_{L^{\frac{2 \gamma}{\gamma - 1}}\left(\cup_i C_i^\eps \right)} + O(\sigmaeps).
\end{align*}
Further, \[
| \nabla(\qk \psi \phi )| \leq | (\nabla \qk )\psi \phi| + | \qk (\nabla \psi \phi + \psi \nabla \phi)| \leq C \left(| \nabla \qk| + \frac{1}{\varepsilon}|\qk| \right),
\]
and hence, \begin{equation} \label{estimate I1}
|I^1| \leq C \sigmaeps^2 \eps^{\frac{\beta}{\gamma}} \left( \| \nabla \qk \|_{L^{\frac{2 \gamma}{\gamma - 1}}\left(\cup_i C_i^\eps \right)}+ \frac{1}{\varepsilon} \| \qk \|_{L^{\frac{2 \gamma}{\gamma - 1}}\left(\cup_i C_i^\eps \right)} \right) + O(\sigmaeps).
\end{equation}

Applying \eqref{eq:qk.C_i} and \eqref{eq:nabla.qk.C_i} from Lemma \ref{lem:qk}  with $ p = (2\gamma)/(\gamma-1)$ yields
$I_1 \to 0$ provided that 
\begin{align} \label{relBeta}
\varepsilon^{3 - \alpha + \frac{\beta}{\gamma} -  2\alpha  +3 \frac{(\alpha - 1)(\gamma-1)}{2\gamma}} \to 0.
\end{align}
Note that since $\alpha > 1$, the term involving $\nabla q_k^\eps$ is always larger than the other term. Hence, $(\ref{relBeta})$ is equivalent to

\[
\beta > \frac{3}{2}(\gamma +1) (\alpha - 1).
\]

%

It remains to show that $|I^2| \to 0$: Since $\psi = 1$ in each $B_i^{\varepsilon/4}$, we have \[
|I^2| = \int_{\omegaeps} \qk (1-\psi) \phi \div(\ueps)\, dx \leq C \sigmaeps \int_{\Omega \setminus \cup_i B_i^{\varepsilon/4}} | \qk| \left|\nabla \frac{\ueps}{\sigmaeps} \right|\, dx \leq C \sigmaeps \| \qk \|_{L^2(\Omega\setminus \cup_i B_i^{\varepsilon/4})}.
\]
By \eqref{eq:qk.annuli} from Lemma \ref{lem:qk},
 \[
|I^2| \leq C \sigmaeps \| \qk \|_{L^2(\Omega\setminus \cup_i B_i^{\varepsilon/4})} \leq C \varepsilon^{\frac{3 - \alpha}{2}} \varepsilon^{\alpha - 2}= C \varepsilon^{\frac{\alpha - 1}{2}} \to 0.
\]
This finishes the proof.
\end{proof}

\section{Homogenization of the unsteady compressible Navier-Stokes equations including a Low-Mach-Number-Limit} \label{sec:Navier.Stokes}

\subsection{A priori estimates}

\begin{lem} \label{aprioriMach3}
Let the assumptions of Theorem \ref{theorem3} hold.
 Then there is a constant $C_T$, which does not depend on $\varepsilon$ such that
\[
\sup_{0 \leq t \leq  T} \int_{\omegaeps} \frac{1}{\eps^\beta} \left| 
\rhoeps^\frac{\gamma}{2}(t) -
\langle \rhoeps \rangle_\e^\frac{\gamma}{2}\right|^{2}\, dx + \sup_{0 \leq t \leq  T} \int_{\omegaeps} \frac{\rhoeps(t) | \ueps^2(t)|}{2}\, dx 
+\left\| \nabla \left( \frac{\utilde}{\sigmaeps} \right)\right\|_{L_T^2(\Omega)}^2
\leq C_T,
\]
and in particular, 
\begin{equation}
\| \rhoeps \|_{L^{\infty}(L^{\gamma}(\omegaeps))} \leq C_T.
\end{equation}
%
\end{lem}


\begin{rem} \label{convergence mean}
The Poincar\'{e} inequality (Lemma~\ref{poincare}) 
implies
 \begin{equation}
\left\|\frac{\utilde}{\sigmaeps^2}\right\|_{L_T^2(\Omega)} \leq c\left\| \nabla \left( \frac{\utilde}{\sigmaeps} \right)\right\|_{L_T^2(\Omega)}\leq C_T.
\end{equation}
\end{rem}
\vspace{0.1 cm}
\begin{proof}[Proof of Lemma~\ref{aprioriMach3}]
Recall from \eqref{energy_eq} that the following energy inequality holds true:
\begin{gather}
\frac{\sigmaeps^2}{\gamma - 1} \frac{1}{\eps^\beta}\int_{\omegaeps} \left( \rhoeps^{\gamma}(T)- \rho_{\varepsilon 0}^{\gamma}\right)\, dx  + \sigmaeps^2 \int_{\omegaeps} \frac{\rhoeps(T)|\ueps(T)|^2}{2}\, dx + \int_0^T \int_{\omegaeps} |\nabla \ueps|^2 + (\div \ueps)^2\, dx\, dt \nonumber \\ 
\leq \sigmaeps^2  \int_{\omegaeps} \frac{\rho_{\varepsilon 0}| u_{\varepsilon 0}|^2}{2}\, dx + \int_0^T \int_{\omegaeps} (\rhoeps f + g) \cdot \ueps\, dx\, dt. \label{energy_eq_3}
\end{gather}

Using that $\langle \rhoeps \rangle_\e = \langle \rho_{\varepsilon 0} \rangle_\e$ by conservation of mass, we rewrite the first term on the left hand side and find by the conservation of the mass and \eqref{eq:gamma2}
\begin{align*}
\int_{\omegaeps} \left( \rhoeps^{\gamma}(T)- \rho_{\varepsilon 0}^{\gamma}\right)\, dx & = \gamma \int_{\omegaeps} \frac{\rhoeps^{\gamma}(T)}{\gamma} - \frac{1}{\gamma} \langle \rhoeps \rangle_\e^{\gamma} - \langle \rhoeps \rangle_\e^{\gamma - 1}  \big(\rhoeps(T)-\langle \rhoeps \rangle_\e \big)\, dx
 \\
& \hspace{3 ex} - \gamma \int_{\omegaeps} \frac{1}{\gamma} \rho_{\varepsilon 0}^{\gamma} - \frac{1}{\gamma} \langle \rho_{\varepsilon 0}  \rangle_\e^{\gamma} - \langle \rho_{\varepsilon 0} \rangle_\e^{\gamma - 1}  \big(\rho_{\varepsilon 0} -\langle \rho_{\varepsilon 0}  \rangle_\e \big)\, dx
\\
& \geq \frac 1 C \int_{\omegaeps} \left| \rhoeps(T)^{\gamma/2} - \langle \rhoeps \rangle_\e^{\gamma/2}\right|^2\, dx - C\int_{\omegaeps} \left| \rho_{\varepsilon 0}^{\gamma/2} - \lrangle{\rho_{\varepsilon 0}}_\eps^{\gamma/2}\right|^2\, dx.
\end{align*}
Thus by assumption \eqref{assumption_rho} we have 
\begin{equation} \label{eq:gamma to gamma/2}
\frac{1}{\eps^\beta} \int_{\omegaeps} \left( \rhoeps^{\gamma}(T)- \rho_{\varepsilon 0}^{\gamma}\right) \geq \frac{1}{C \eps^\beta} \int_{\omegaeps} \left| \rho_{\varepsilon}^\frac{\gamma}{2}(T) - \langle \rhoeps \rangle_\e ^\frac{\gamma}{2}\right|^{2} - C.
\end{equation}
We next estimate the right-hand side of the energy inequality.
In case $\gamma\geq 2$ we find for $b$ satisfying $\frac{1}{b}= \frac{1}{2}-\frac{1}{\gamma}$, using H\"older's and Young's inequalities and Lemma \ref{poincare}, that
\begin{align*}
\int_{\omegaeps} |(\rhoeps f + g) \cdot u_{\varepsilon}|\, dx &
\leq  \int_{\omegaeps} \left|\sqrt{2C\sigmaeps^2}(\dpr{\rhoeps}_\e)f + g) \cdot \sqrt{\frac{1}{2C \sigmaeps^2}}u_{\varepsilon}\right|+ \abs{(\rhoeps -\dpr{\rhoeps}_\e)}\;\abs{f\ueps}\, dx
\\
& \leq C \sigmaeps^2 \left(m_0^2\|f\|_{L^{2}(\omegaeps)}^2 + \|g\|_{L^{2}(\omegaeps)}^2\right) + \frac{1}{C \sigmaeps^2} \| \ueps\|_{L^2(\omegaeps)}^2
\\
&\quad + \sigmaeps\norm{\rhoeps -\dpr{\rhoeps}_\e}_{L^\gamma(\omegaeps)}\norm{f}_{L^b(\omegaeps)}\normB{\frac{\ueps}{\sigmaeps}}
_{L^2(\omegaeps)}
\\
& \leq C \sigmaeps^2 \left(m_0^2\|f\|_{L^{2}(\omegaeps)}^2 + \norm{f}_{L^b(\omegaeps)}^{b}+\|g\|_{L^{2}(\omegaeps)}^2\right) + \frac{1}{2} \| \nabla \ueps\|_{L^2(\omegaeps)}^2
\\
&\quad + C\sigmaeps^2\norm{\rhoeps^\frac{\gamma}{2} -\dpr{\rhoeps}_\e^\frac{\gamma}{2}}_{L^2(\omegaeps)}^2,
\end{align*}
if $b < \infty$, i.e $\gamma \neq 2$. For $\gamma = 2$, $\|f\|_\infty$ appears as a factor in front of the last term on the right-hand side and does not appear in the first term.

In order to estimate the right-hand sides in case $\gamma<2$ and $f\neq 0$ we have to proceed differently. Indeed, using $\rho_\eps = (\sqrt{\rho_\eps} - \sqrt{\lrangle{\rho_\eps}_\eps})(\sqrt{\rho_\eps} +  \sqrt{\lrangle{\rho_\eps}_\eps}) + \lrangle{\rho_\eps}_\eps$ and Young's inequality, we find
\begin{align*}
|\rhoeps f\cdot u_{\varepsilon}|
&\leq \norm{f}_\infty\abs{\rhoeps\ueps}
\leq \norm{f}_\infty \big( \abs{\sqrt{\rhoeps}}\abs{\ueps}\abs{\rhoeps^\frac{\gamma}{2}-\dpr{\rhoeps}_\e^\frac{\gamma}{2}}^\frac{1}{\gamma}
+\abs{\sqrt{\dpr{\rhoeps}_\e}}\abs{\ueps}\abs{\rhoeps^\frac{\gamma}{2}-\dpr{\rhoeps}_\e^\frac{\gamma}{2}}^\frac{1}{\gamma}+\dpr{\rhoeps}_\e\abs{\ueps}\big)
\\
&\leq {\sigmaeps^2}\rhoeps \abs{\ueps}^2+\frac{C\abs{\rhoeps^\frac{\gamma}{2}-\dpr{\rhoeps}_\e^\frac{\gamma}{2}}^\frac{2}{\gamma}}{\sigmaeps^2} + \frac{1}{C}\left| \frac{\ueps}{\sigmaeps}\right|^2+2C\sigmaeps^2m_0\abs{\rhoeps^\frac{\gamma}{2}-\dpr{\rhoeps}_\e^\frac{\gamma}{2}}^\frac{2}{\gamma} +Cm_0\sigmaeps^2
\\
&\leq {\sigmaeps^2}\rhoeps \abs{\ueps}^2+\frac{C\sigmaeps^2\abs{\rhoeps^\frac{\gamma}{2}-\dpr{\rhoeps}_\e^\frac{\gamma}{2}}^2}{\sigmaeps^{4 \gamma}} + \frac{1}{C_p}\left| \frac{\ueps}{\sigmaeps}\right|^2+2C\sigmaeps^2m_0\abs{\rhoeps^\frac{\gamma}{2}-\dpr{\rhoeps}_\e^\frac{\gamma}{2}}^2 +C\sigmaeps^2,
\end{align*}
where $C_p$ is the Poincaré constant from Lemma \ref{poincare} and the constant $C$ depends on $f$ and $m_0$.
Here, we need to use the assumption that $2 \gamma(3-\alpha) < \beta$ which implies that $\frac{\eps^\beta}{\sigma_\eps^{4\gamma}}\to 0$.
Thus, we find in all assumed cases that for $\eps$ sufficiently small
\begin{align*}
&\int_0^T \int_{\omegaeps} |(\rhoeps f + g) \cdot u_{\varepsilon}|\, dx\, dt \\
& \leq C \sigmaeps^2 + \frac{1}{2} \| \nabla \ueps\|_{L_T^2(\omegaeps)}^2
 +C \int_0^T{\sigmaeps^2}\norm{\rhoeps\abs{\ueps}^2}_{L^1(\omegaeps)}+ \frac{\sigmaeps^2}{\eps^\beta(\gamma-1)}\norm{\rhoeps^\frac{\gamma}{2} -\dpr{\rhoeps}_\e^\frac{\gamma}{2}}_{L^2(\omegaeps)}^2\, dt.
\end{align*}

Combining this estimate with \eqref{eq:gamma to gamma/2}, the energy inequality \eqref{energy_eq_3}{} yields
\begin{align*}
&\frac{1}{\gamma - 1}\sigmaeps^2 \int_{\omegaeps} \frac{1}{\eps^\beta}\absB{\rhoeps^\frac{\gamma}{2}(T)-\dpr{\rhoeps}_\e^\frac{\gamma}{2}}^2+ \sigmaeps^2\frac{\rhoeps (T) |\ueps (T)|^2}{2}\, dx + \frac12 \int_0^T \int_{\omegaeps} |\nabla \ueps|^2\, dx\, dt.
\\
& \quad \leq C  \sigmaeps^2 \left(1 + \frac{1}{\eps^\beta}\int_0^T \int_{\omegaeps} \absB{\rhoeps^\frac{\gamma}{2}-\dpr{\rhoeps}_\e^\frac{\gamma}{2}}^2\, dx\, dt\right)+  \sigmaeps^2 \int_0^T \int_{\omegaeps} {\rhoeps | \ueps|^2}\, dx\, dt + \sigmaeps^2 \int_{\omegaeps} \left( \frac{\rho_{\varepsilon 0} |u_{\varepsilon 0}|^2}{2} \right)\, dx
\end{align*}

Dividing by $\sigmaeps^2$, we conclude by Gronwall's Theorem that for all $T>0$
\begin{align} \label{l infty estimate}
\sup_{0 \leq t \leq  T} \int_{\omegaeps} \frac{1}{\eps^\beta} \absB{\rhoeps^\frac{\gamma}{2}-\dpr{\rhoeps}_\e^\frac{\gamma}{2}}^2\, dx+ \sup_{0 \leq t \leq  T} \int_{\omegaeps} \frac{\rhoeps(t) | \ueps(t)|^2}{2}\, dx+ \int_0^T \int_{\omegaeps} \left|  \frac{\nabla\ueps}{\sigmaeps} \right|^2 \, dx\, dt\leq C_T.
\end{align}
For $\rhoeps$, this in particular means \begin{equation}
\| \rhoeps\|_{L^{\infty}_T(L^{\gamma}(\omegaeps))} \leq \|\rhoeps^\frac{\gamma}{2}-\dpr{\rhoeps}_\e^\frac{\gamma}{2}\|^{\frac 2 \gamma}_{L^{\infty}_T(L^{2}(\omegaeps))} + C_T \langle \rhoeps \rangle_\e \leq C_T,
\end{equation}
which concludes the proof.
\end{proof}

In order to pass to the limit as it was done in the last section, we are left to estimate the pressure, which we will now decompose as stated in Theorem \ref{theorem3}.
To this end, we first define a modified Bogovski\u{\i} operator that also acts on functions 
that are not mean free by first removing the average. More precisely, we define
%
\begin{align}
B^*_\e:L^a(\Omega)\to W^{1,a}_0(\Omega_\e),\quad  B^*_\e(\phi)=\Bcal_\e((\phi-\dpr{\phi}_\e)\chi_{\Omega_\e}),
\end{align}
where $\chi_{\Omega_\e}$ is the characteristic function of the set $\omegaeps$.

First we find for $\phi\in C^1_0((0,T)\times\Omega))$
\begin{align*}
\langle \tilde p_\e,\phi\rangle_{\mathcal{D}', \mathcal{D}((0,T)\times \Omega)}&
= 
\int_0^T \int_{\omegaeps}\left(\frac{\rho_{\varepsilon}^{\gamma}}{\varepsilon^{\beta}} - \langle \frac{\rho_{\varepsilon}^{\gamma}}{\varepsilon^{\beta}}\rangle_{\varepsilon}\right) \phi \, dx \, dt = \int_0^T \int_{\omegaeps}\frac{\rho_{\varepsilon}^{\gamma}}{\varepsilon^{\beta}}(\phi - \langle \phi \rangle_{\varepsilon}) \, dx \, dt\\
& = \int_0^T \int_{\omegaeps}\frac{\rho_{\varepsilon}^{\gamma}}{\varepsilon^{\beta}}\div( \mathcal{B}_{\varepsilon}(\phi - \langle \phi \rangle_{\varepsilon})) \, dx \, dt 
 = \int_0^T \int_{\omegaeps}\frac{\rho_{\varepsilon}^{\gamma}}{\varepsilon^{\beta}}\div( \mathcal{B}_{\varepsilon}^{*}(\phi)) \, dx \, dt\\
& =-\sigmaeps^2 \int_{0}^T \int_{\omegaeps} \rhoeps \ueps \partial_{t} B^*_\e(\phi) \, dx\, dt - \int_0^T \int_{\omegaeps} \rhoeps \ueps \otimes \ueps : \nabla B^*_\e(\phi)\, dx\, dt\\
& \hspace{2 ex} + \int_0^T \int_{\omegaeps} \nabla \ueps : \nabla B^*_\e(\phi)+ \div(\ueps) \div(B^*_\e(\phi))- (\rhoeps f + g)B^*_\e(\phi)\, dx\, dt
\end{align*}
and split the pressure into two parts  $\tilde p_\e=p_{\e,1}+p_{\e,2}$, where 
\begin{align}
\label{eq:p1}
\langle p_{\e,1},\phi\rangle_{\mathcal{D}', \mathcal{D}((0,T)\times \Omega)}:=-\sigmaeps^2 \int_{0}^T \int_{\omegaeps} \rhoeps \ueps \partial_{t} B^*_\e(\phi) \, dx\, dt - \int_0^T \int_{\omegaeps} \rhoeps \ueps \otimes \ueps : \nabla B^*_\e(\phi)\, dx\, dt
\end{align}
and
\begin{align}
\label{eq:p2}
\langle p_{\e,2},\phi\rangle_{\mathcal{D}', \mathcal{D}((0,T)\times \Omega)}&:= \int_0^T \int_{\omegaeps} \nabla \ueps \cdot\nabla B^*_\e(\phi)+ \div(\ueps) \div(B^*_\e(\phi)) -  (\rhoeps f + g)B^*_\e(\phi)\, dx\, dt.
\end{align}
Since $p_{\e,2}$ basically coincides with $p_\e$ as defined in \eqref{eq:defp2} (for a.e.\ $t\in [0,T]$), we find that it is bounded in $L^2([0,T]\times(\Omega))$ by similar arguments as in the previous section. This will be shown rigorously in the next lemma. 
Furthermore, in Lemma~\ref{lem:p1} we will show that $p_{\varepsilon,1}$ converges to $0$ in some appropriate topology.
\begin{lem}
\label{lem:p2}
We find that $\norm{p_{\e,2}}_{L^2([0,T]\times \Omega)}\leq C$, with a constant independent of $\e$.
\end{lem}
\begin{proof}
By H\"older's inequality and Sobolev embedding, the a priori estimates and Theorem \ref{thm_bog_eps}, we find that in case $\gamma\geq 2$
\begin{align*}
&\langle p_{\e,2},\phi\rangle_{\mathcal{D}', \mathcal{D}((0,T)\times \Omega)}
\\
&\quad \leq  \int_0^T\norm{\nabla \ueps}_{L^2(\omegaeps)}\norm{\nabla B_\e^*(\phi)}_{L^2(\omegaeps)}+(\norm{\rhoeps}_{L^{\gamma}(\omegaeps)}\norm{f}_{L^b(\omegaeps)}+\norm{g}_{L^2(\omegaeps)})\norm{B_\e^*(\phi)}_{L^2(\omegaeps)}\, dt
\\
&\quad \leq  C \norm{\phi}_{L^2([0,T]\times \Omega)} .
\end{align*}
In case $\gamma\in (\frac32,2)$ and $f\neq 0$ we have to estimate the term involving $f$ differently. Observe that  there is a constant $C$ only depending on $\gamma$ such that for all $a,b >0$
\[
a\leq C\abs{a^\frac{\gamma}{2}-b^\frac{\gamma}{2}}^\frac{2}{\gamma}+Cb.
\]

Thus,
\begin{align*}
\int_0^T\int_{\omegaeps} \rhoeps f\cdot B_\e^*(\phi)\, dx\, dt
\leq & C\int_0^T\int_{\omegaeps} \Big(\abs{\rho_\eps^\frac{\gamma}{2}-\lrangle{\rho_\eps}^\frac{\gamma}{2}}^\frac{2}{\gamma}+\lrangle{\rho_\eps}_\eps\Big)
| f\cdot B_\e^*(\phi)|\, dx\, dt
\\
 \leq &  C\int_0^T\norm{\rho_\eps^\frac{\gamma}{2}-\lrangle{\rho_\eps}_\eps^\frac{\gamma}{2}}_{L^{2}(\omegaeps)}^{\frac 2 \gamma}
\norm{f}_{L^\infty(\omegaeps)}\norm{B_\e^*(\phi)}_{L^{\gamma'}(\omegaeps)}\, dt \\
&+ C\int_0^T(\dpr{\rho_\eps}_\e\norm{f}_{L^2(\omegaeps)}
\norm{B_\e^*(\phi)}_{L^2(\omegaeps)} \, dt
\\
 \leq & C \eps^\frac{\beta}{\gamma} \int_0^T\norm{B_\e^*(\phi)}_{L^{\gamma'}(\omegaeps)} +  C\norm{\phi}_{L^2([0,T]\times \Omega)}\, dt
\\
 \leq & C\norm{\phi}_{L^2([0,T]\times \Omega)},
\end{align*}
where we used that  Sobolev embedding and Theorem~\ref{thm_bog_eps} imply
\[
\norm{B_\e^*(\phi)}_{L^{\gamma'}(\omegaeps)}\leq c\norm{\nabla B_\e^*(\phi)}_{L^\frac{3\gamma'}{3+\gamma'}(\omegaeps)}\leq c\norm{\nabla B_\e^*(\phi)}_{L^{2}(\omegaeps)}\leq c\sigmaeps^{-1}\norm{\phi}_{L^2(\Omega)}\leq C\eps^{-\frac{\beta}\gamma}\norm{\phi}_{L^2(\Omega)},
\]
whenever $\beta \geq \frac{\gamma}{2} (3 - \alpha)$.
%
%
\end{proof}

\subsection{Convergence Proof}
The most significant difference to the last section is the appearance of the new pressure term $p_{1,\e}$. The following lemma shows why the terms related to that quantity vanish:
\begin{lem}
\label{lem:p1}
There exists $q\in (1,2]$ such that
$p_{\e,1} \to 0$, strongly  in $H^{-1}(0,T;L^q(\Omega))$, for $p_{1,\e}$ defined in \eqref{eq:p1}. Moreover, 
\begin{align} \label{eq:p_1.w_k.to.0}
\langle p_{\e,1},\div(\wk \phi)\rangle_{\mathcal{D}', \mathcal{D}((0,T)\times \Omega)}\to 0
\end{align}
for all $\phi\in C^\infty_0((0,T)\times\Omega)$.
\end{lem}
\begin{proof}
We begin with the proof of \eqref{eq:p_1.w_k.to.0}. Observe, that $\int_{\omegaeps}\div(\wk \phi)\, dx=0$ and $\div(\wk \phi) = \wk \cdot \nabla \phi$. Hence,
\begin{align*}
\langle p_{\e,1},\div(\wk \phi)\rangle_{\mathcal{D}', \mathcal{D}((0,T)\times \Omega)}&=-\sigmaeps^2 \int_{0}^T \int_{\omegaeps} \rhoeps \ueps \Bcal_\e(\wk\cdot \nabla \partial_t\phi) \, dx\, dt
\\
&\quad - \int_0^T \int_{\omegaeps} \rhoeps \ueps \otimes \ueps : \nabla \Bcal_\e(\wk\cdot \nabla \phi)\, dx\, dt \\
&=: B_1+B_2. 
\end{align*}
With $\frac 1 3 + \frac 1 \gamma + \frac 1 a = 1$, we estimate
\begin{align*}
\abs{B_2}&\leq C \int_0^T \int_{\Omega_\eps} \Big(\abs{\rho_\eps^\frac{\gamma}{2}-\lrangle{\rho_\eps}_\eps^\frac{\gamma}{2}}^\frac{2}{\gamma}+ \lrangle{\rho_\eps}_\eps\Big)\abs{u_\eps}^2\abs{ \nabla \Bcal_\e(\wk\cdot \nabla \phi)}\, dx\, dt
\\
&\leq C \sigmaeps^2 \int_0^T \norm{\rho_\eps^\frac{\gamma}{2}-\lrangle{\rho_\eps}_\eps^\frac{\gamma}{2}}^{\frac 2 \gamma}_{L^2(\omegaeps)}\normB{\frac{u_\eps}{\sigmaeps}}_{L^6(\omegaeps)}^2\norm{\nabla \Bcal_\e(\wk\cdot \nabla \phi)}_{L^a(\omegaeps)}\, dt
\\
&\quad  + C m_0\sigmaeps^2\int_0^T\normB{\frac{u_\eps}{\sigmaeps}}_{L^4(\omegaeps)}^2\norm{\nabla B_\e(\wk\cdot \nabla \phi)}_{L^2(\omegaeps)}\, dt.
\end{align*}
Using Theorem~\ref{thm_bog_eps} and the fact that $\wk$ is uniformly bounded, we get
\[
\norm{\nabla B_\e(\wk\cdot \nabla \phi)}_{L^a(\omegaeps)}\leq C \e^\frac{(3-a)\alpha-3}{a}\norm{\wk}_{L^\infty(\Omega)}\norm{\nabla \phi}_{L^a(\Omega)}
\leq C\e^\frac{(3-a)\alpha-3}{a}.
\]
Hence we find by the definition of $\beta$ that
\begin{equation}
|B_2| \leq  C\e^{(3-\alpha)}\e^\frac{(3-a)\alpha-3}{a}\eps^\frac{\beta}{\gamma} +C\sigmaeps\to 0
\end{equation}
with $\e\to 0$, due to the assumption on $\beta$.
For $B_1$ let $b < 6$ satisfy $\frac{1}{b}+\frac{1}{6}+\frac{1}{\gamma}=1$. Then, by Sobolev inequality,
\begin{align*}
|B_1| &\leq C \sigmaeps^3 \int_0^T \norm{\rho_\eps}_{L^\gamma(\omegaeps)}
\normB{\frac{u_\eps}{\sigmaeps}}_{L^6(\omegaeps)}
\norm{\Bcal_\e (\wk \nabla \partial_t\phi)}_{L^b(\omegaeps)}\, dt\\
&\leq C \sigmaeps^3 \int_0^T \norm{ \nabla \Bcal_\e (\wk \nabla \partial_t\phi)}_{L^2(\omegaeps)} \, dt \leq C \sigma_\eps^2  \int_0^T  \|\nabla \partial_t\phi\|_{L^2(\Omega)}\, dt \to 0
\end{align*}
with $\e\to 0$. This proves \eqref{eq:p_1.w_k.to.0}.


 The first assertion of the lemma, $p_{\eps,1} \to 0$, follows by very similar arguments. 
More precisely, for $\phi \in \Dcal(\Omega_T)$,
\begin{align}
\langle p_{\e,1},\phi\rangle&=-\sigmaeps^2 \int_{0}^T \int_{\omegaeps} \rhoeps \ueps B^*_\e(\partial_{t} \phi)  + \rhoeps \ueps \otimes \ueps : \nabla B^*_\e:(\phi)\, dx\, dt =: C_1 + C_2.
\\
\end{align}
In analogy of the above we find
(using the linearity of all involved operators, Sobolev embedding and Theorem~\ref{thm_bog_eps})
\begin{align}
|C_1|\leq C \sigmaeps^2 \| \partial_t \phi \|_{L^2(\Omega)},
\end{align}
and 
\begin{align}
|C_2| \leq C \sigmaeps^2 \left( \varepsilon^ {\frac{(3 - a)\alpha - 3}{a}} \varepsilon^{\frac{\beta}{\gamma}} \| \phi \|_{L_T^\infty(L^a(\Omega))} + \frac{1}{\sigmaeps} \| \phi \|_{L_T^\infty(L^2(\Omega)} \right). 
\end{align}
Using the embedding $H^1_T \subset L^\infty_T$, these estimates imply the assertion.
\end{proof}

\begin{proof}[Proof of Theorem \ref{theorem3}]

\emph{Step 1:} Extracting weakly convergent subsequences.

We use the a priori estimates and Lemma \ref{lem:p1} to find  $u \in L^2_T(\Omega)$  and $p \in L^2_T(L^2_0(\Omega))$ and  subsequences (not relabeled) such that

\begin{alignat*}{2}
\rhotilde &\to  \rho_0  \quad && \text{ strongly in }L^{\infty}_T( L^{ \gamma}(\Omega)) 
,\\
\frac{\utilde}{\sigmaeps^2} &\weak u  \quad && \text{ weakly in } L_T^{2}(\Omega),\\
p_{\e,1} &\to 0\quad && \text{ strongly in }H^{-1}(0,T;L^q(\Omega)),
\\
p_{\e,2}&\weakto p\quad && \text{ weakly in }L^2([0,T]\times\Omega)).
\end{alignat*}
Again, if we prove that $(u,p)$ solves \eqref{limit_system1}, convergence of the full sequence follows.

\medskip

\emph{Step 2:} Passage to the limit in the continuity equation.

Recall from \eqref{eq:continuity.everywhere.assumption} that the continuity equation holds for $(\tilde \rho_\eps,\tilde u_\eps)$ in $\R^3 \times (0,T)$.
We test the equation with $\phi \in \Dcal(\R^3 \times (0,T))$. For $\gamma \geq 2$, we can directly pass to the limit to deduce
\begin{align}
	0 &= \int_0^T \int_\Omega \tilde \rho_\eps \partial_t \phi+\tilde \rho_\eps \frac{\tilde u_\eps}{\sigma_\eps^2} \cdot \nabla \phi \, dx\, dt \to \int_0^T \int_\Omega \rho_0 \left(\partial_t \phi + u \cdot \nabla \phi\right)\, dx\, dt = \int_0^T \int_\Omega \rho_0 u \cdot \nabla \phi\, dx\, dt.
\end{align} 
Thus ,
\begin{equation} \label{eq:u.dviergencefree}
\begin{cases}
	\div u = 0 &\quad \text{in } \Omega \\
	u \cdot n = 0 &\quad \text{on } \partial \Omega.
\end{cases}
\end{equation}

For $\gamma < 2$, we again have to proceed differently.
We rewrite the critical term
\begin{align}
	\int_0^T \int_\Omega\tilde \rho_\eps \frac{\tilde u_\eps}{\sigma_\eps^2} \cdot \nabla \phi \, dx\, dt = 
	\int_0^T \int_{\Omega_\eps} (\rho_\eps - \lrangle{\rho_\eps}_\eps) \frac{ u_\eps}{\sigma_\eps^2} \cdot \nabla \phi \, dx\, dt + \int_0^T \int_\Omega \langle \rho_\eps\rangle_\eps \frac{\tilde u_\eps}{\sigma_\eps^2} \cdot \nabla \phi \, dx\, dt.
\end{align}
Then, we may pass to the limit in the second term on the right-hand side and it remains to prove that the first term on the right-hand side vanishes in the limit $\eps \to 0$. 
To this end, we use $\rho_\eps - \lrangle{\rho_\eps}_\eps = (\sqrt{\rho_\eps} - \sqrt{\lrangle{\rho_\eps}_\eps}) (\sqrt{\rho_\eps} + \sqrt{\lrangle{\rho_\eps}_\eps})$ to estimate 
\begin{align*}
	\left\|(\rho_\eps - \lrangle{\rho_\eps}_\eps) \frac{ u_\eps}{\sigma_\eps^2}\right\|_{L^1_T(\Omega_\eps)} & \leq 
	\int_0 ^T \|\rho_\eps^{\frac \gamma 2} - \lrangle{\rho_\eps}_\eps^{\frac \gamma 2}\|_{L^2(\Omega_\eps)}^{\frac 1 \gamma}\left( \sigma_\eps^{-2}\|\rho_\eps u_\eps^2\|_{L^1(\Omega_\eps)}^{\frac 1 2} +  \sqrt{\lrangle{\rho_\eps}_\eps}  \left\|\frac{\ueps}{\sigmaeps^2} \right\|_{L^2(\Omega_\eps)} \right) \\
	& \leq C \eps^{\frac \beta {2 \gamma}} \eps^{-3 + \alpha} \to 0,
\end{align*}
by the strengthened assumptions on $\beta$ for $\gamma < 2$. 
This shows that \eqref{eq:u.dviergencefree} also holds for $\gamma < 2$.

\medskip

\emph{Step 3:} Passage to the limit in the momentum equation.

Let $\phi \in \mathcal{D}((0,T)\times \Omega)$. Using $\wk \phi$ as test function in the first equation of \eqref{goal}, rewriting the terms in the same way as before and adding the term $ - \sigma_{\varepsilon}^2 \int_0^T \int_{\Omega}  \qk \div \left(\frac{\utilde}{\sigma_{\varepsilon}^2} \phi  \right)$ to both sides gives

\begin{align*}
&\underbrace{\sigma_{\varepsilon}^2 \int_0^T\int_{\Omega} \nabla \wk :\nabla\left(\frac{\utilde}{\sigma_{\varepsilon}^2} \phi \right)-  \qk \div \left(\frac{\utilde}{\sigma_{\varepsilon}^2} \phi\, dx\, dt \right)}_{I_1} 
\\
&- \underbrace{\int_0^T \int_{\Omega} (\rhotilde \utilde\otimes \utilde): \nabla(\wk \phi)\, dx\, dt}_{I_2} + \underbrace{\int_0^T \int_{\Omega} \div(\utilde) \div(\wk \phi) \, dx\, dt}_{I_3} 
\\
&+ \underbrace{\int_0^T\int_{\Omega} \nabla \utilde :  (\wk \otimes \nabla \phi) -  \nabla \wk:(\utilde \otimes \nabla \phi)\, dx\, dt}_{I_4} - \underbrace{\sigmaeps^2 \int_0^T\int_{\Omega} \rhotilde \utilde\cdot \wk \partial_t\phi\, dx\, dt}_{I_5}
\\
&= \underbrace{\frac{1}{\eps^\beta}\int_0^T \int_{\Omega} (\rhoeps^\gamma-\dpr{\rhoeps^\gamma}_\e) \div(\wk \phi)\, dx\, dt}_{I_6} +\underbrace{ \int_0^T \int_{\Omega} (\rhoeps f + g)(\wk \phi)\, dx\, dt}_{I_7} -\underbrace{ \int_0^T \int_{\Omega} \qk \div \left(\utilde \phi \right)\, dx\, dt}_{I_8}.
\end{align*}

Analogously to the convergence proofs in Theorem~\ref{theorem1} we find
\begin{gather*}
I_1 \to \int_0^T \int_{\Omega} \phi \Rcal e_k \cdot u\, dx\, dt,
\end{gather*}
$ I_3, \, I_4 \to 0$ and 
\[
I_7 \to \int_0^T \int_{\Omega} (\rho_0 f + g)\cdot e_k \phi\, dx\, dt
\]
as $\varepsilon \to 0$. Furthermore, we choose $a = \frac{3\gamma}{2\gamma - 3}$ such that $\frac{1}{\gamma}+\frac{1}{3}+\frac{1}{a}=1$ and find by triangular and H\"older's inequality
\begin{align*}
\abs{I_2}&\leq c\int_0^T \int_{\Omega}\Big( \abs{\rhotilde^\frac\gamma2-\lrangle{\rhotilde}^\frac\gamma2}^\frac2\gamma+\lrangle{\rhotilde}\Big)\abs{\utilde}^2\abs{\nabla(\wk \phi)}\, dx\, dt
\\
&\leq c\sigmaeps^2 \int_0^T \norm{\rhotilde^\frac\gamma2-\lrangle{\rhotilde}^\frac\gamma2}^{\frac 2 \gamma}_{L^2(\omegaeps)}\normB{\frac{\utilde}{\sigmaeps}}_{L^6(\omegaeps)}^2\norm{\nabla(\wk \phi)}_{L^a(\Omega)}\, dt \\
&+m_0\sigmaeps^2\int_0^T\normB{\frac{\utilde}{\sigmaeps}}_{L^4(\omegaeps)}^2\norm{\nabla(\wk \phi)}_{L^2(\Omega)}\, dt
\\
&=: I_2^1+I_2^2 .
\end{align*}
The term $I_2^2$ can be estimated using the a priori estimates and hypothesis $(H3^{*})$ in Theorem \ref{hypotheses}:
$
I_2^2 \leq C\sigmaeps  \to 0.
$
On the term $I_2^1$ we use the bounds on $\beta$ and Lemma~\ref{aprioriMach3} together with the uniform $L^\infty$-bounds on $w_k^\eps$ and the $L^a$-bounds on $\nabla w_k^\eps$ from Lemma \ref{lem:qk} to find
\[
I_2^1\leq C \eps^{3 - \alpha} \varepsilon^\frac{\beta}{\gamma} \varepsilon^{- \alpha + 3(\alpha-1)\frac{2\gamma -3}{3 \gamma}} \to 0.
\]
%


Further,
\begin{align*}
|I_5| &\leq C \sigmaeps^3 \int_0^T \norm{\rhotilde}_{L^\gamma(\omegaeps)}\normB{\frac{\utilde}{\sigmaeps}}_{L^6(\omegaeps)}\norm{\wk }_{L^a(\Omega)}\norm{\partial_t\phi}_{L^6(\omegaeps)}
\to 0
\end{align*}
as $\varepsilon \to 0$. 

Since $\wk \to e_k$ strongly in $L^2(\Omega)$,  Lemma~\ref{lem:p2} and Lemma~\ref{lem:p1}
imply
\[
I_6 = \int_0^T \int_{\Omega} (p_{\varepsilon, 1}+p_{\varepsilon,2}) \wk \cdot \nabla \phi \, dx\, dt\to \int_0^T \int_{\Omega} p \, e_k \cdot \nabla \phi\, dx\, dt.
\]

We postpone the proof that $I_8 \to 0$ to the end of the proof.
Accepting this convergence yields
\begin{equation}
 \int_0^T \int_{\Omega} \left(\rho_0f + g-\nabla p \right) e_k \phi\, dx\, dt= \int_0^T\int_{\Omega} \phi \Rcal e_k \cdot u \, dx\, dt
\end{equation}
for all $\phi \in \mathcal{D}((0,T)\times\Omega)$, which means \[
\Rcal u = \rho_0 f + g - \nabla p.
\]
Combining this with \eqref{eq:u.dviergencefree} yields the limit equation \eqref{limit_system1}.


\medskip

\emph{Step 4:} Estimate for $I_8$.

This step corresponds to step 4 in the proof of Theorem \ref{theorem1}.
We split $I_8$ into the analogous terms $I^1, I^2, I^3$ the only difference being that these terms contain also integrals over time.
Since the spacial a priori estimates for $u_\eps$ are identical to the ones in the stationary case, $I^2 \to 0$ and $I^3 \to 0$ is established as before.

The treatment of $I^1$ is slightly different. We use \eqref{eq:continuity.everywhere.assumption} to obtain, analogously as in the stationary case, 
\begin{align*}
I^1 & =  \int_0^T \int_{\omegaeps} \nabla (\qk \psi \phi) \rhoeps \ueps \, dx\,dt
+  \sigmaeps^2 \int_0^T \int_{\omegaeps} \qk \psi \partial_t\phi  \rhoeps 
\, dx\,dt- \int_0^T \int_{\omegaeps} \nabla (\qk \psi \phi) \langle \rhoeps \rangle \ueps \, dx\,dt 
\\
&\quad + \langle \rhoeps \rangle \sigmaeps^2 \int_0^T \int_{\omegaeps} \qk \psi  \frac{\ueps}{\sigmaeps^2} \nabla \phi\, dx\,dt
= \sigmaeps^2 \int_0^T \int_{\omegaeps} \nabla (\qk \psi \phi) (\rhoeps -\langle \rhoeps \rangle_\eps) \,  \frac{\ueps}{\sigmaeps^2}\, dx\,dt + O(\sigmaeps).
\end{align*}
Thus, for $\gamma \geq  2$
\begin{align}
	\limsup_{\eps \to 0} |I^1| &\leq \limsup_{\eps \to 0} \sigmaeps^2 \|\rhoeps^{\frac \gamma 2} -\langle \rhoeps \rangle_\eps^{\frac \gamma 2}\|_{L^\infty_T(L^2(\Omega))}^{\frac 2 \gamma} \left\|\frac{\ueps}{\sigmaeps^2} \right\|_{L^2(\Omega_T)}  \|\nabla (\qk \psi \phi)\|_{L^2_T(L^{\frac{2 \gamma}{\gamma - 2}}(\Omega))} \\
	&\leq C \limsup_{\eps \to 0} \eps^{\frac \beta \gamma} \eps^{3 - \alpha -2 \alpha + 3 (\alpha - 1) \frac{\gamma - 2}{2 \gamma}}
\end{align}
due to Lemma \ref{lem:qk} since $\supp \psi \subset \cup_i C_i^\eps$ .
Thus, $I^1 \to 0$ by the assumptions on $\beta$.

Again, for $\gamma < 2$, we need to proceed slightly different. 
We use $\rho_\eps - \lrangle{\rho_\eps} = (\sqrt{\rho_\eps} - \sqrt{\lrangle{\rho_\eps}}) (\sqrt{\rho_\eps} + \sqrt{\lrangle{\rho_\eps}})$ to estimate 
\begin{align}
	 |I^1| \leq \int_0 ^T \|\rho_\eps^{\frac \gamma 2} - \lrangle{\rho_\eps}_\eps^{\frac \gamma 2}\|_{L^2(\Omega_\eps)}^{\frac 1 \gamma}\left( \|\rho_\eps u_\eps^2\|_{L^1(\Omega_\eps)}^{\frac 1 2} + \sigma_\eps^2 \sqrt{\lrangle{\rho_\eps}_\eps}  \left\|\frac{\ueps}{\sigmaeps^2} \right\|_{L^2(\Omega_\eps)} \right)  \|\nabla (\qk \psi \phi)\|_{L^{\frac{2 \gamma}{\gamma - 1}}(\Omega)} + O(\sigma_\eps^2).
\end{align}
Thus,
\begin{align}
	\limsup_{\eps \to 0} |I^1| \leq \limsup_{\eps \to 0} \eps^{\frac{\beta}{2\gamma}} \eps^{- 2\alpha + 3 (\alpha - 1) \frac{\gamma - 1}{2 \gamma}} = 0
\end{align}
by the strengthened assumption on $\beta$ for $\gamma < 2$. 
This finishes the proof.

\end{proof}

\section*{Acknowledgements}
The authors thank the support of the Hausdorff Center of Mathematics at the University of Bonn.
R.M. H\"ofer also acknowledges the support of the Deutsche Forschungsgemeinschaft (DFG, German Research Foundation) through the collaborative research center ``The Mathematics of Emerging Effects'' (CRC 1060, Projekt-ID 211504053). S. Schwarzacher thanks the support of the Primus research programme PRIMUS/19/SCI/01 and the University Centre UNCE/SCI/023 of Charles University. Moreover he thanks for the support of the program GJ19-11707Y of the Czech national grant agency (GA\v{C}R).

\section*{Conflict of interest}
The authors declare that there is no conflict of interest regarding the publication of this article.

\printbibliography

\end{document}